\let\ovec\vec
\documentclass[smallextended]{svjour3}       
\let\vec\ovec

\smartqed  

\journalname{Journal}
\usepackage{epsfig}
\usepackage{times}
\newcommand{\epc}{\hspace{1pc}}
\newcommand{\thalf}{{\textstyle{\frac{1}{2}}}}
\newcommand{\wh}{\widehat}
\newcommand{\wt}{\widetilde}

\SetSymbolFont{letters}{normal}{OML}{txmi}{m}{it}
\SetSymbolFont{operators}{normal}{OT1}{txr}{m}{n}
\SetSymbolFont{largesymbols}{normal}{OMX}{txex}{m}{n}
\SetMathAlphabet{\mathrm}{normal}{OT1}{txr}{m}{n}

\usepackage[b5paper,left=1.8cm,right=1.8cm,top=2.95cm,bottom=2.95cm]{geometry}
\usepackage{amsmath,amssymb}
\usepackage{mathabx,float}
\usepackage{graphicx}
\usepackage{color}
\usepackage[colorlinks,anchorcolor = black,linkcolor = black,citecolor = blue]{hyperref}
\usepackage{url}
\usepackage{enumitem}
\usepackage{tikz}
\usetikzlibrary{shapes,arrows,matrix,decorations.markings}

\definecolor{blue}{rgb}{0,0,0.9}
\definecolor{red}{rgb}{0.9,0,0}
\definecolor{green}{rgb}{0,0.9,0}

\renewcommand\zeta{f}
\renewcommand\nu{g}

\begin{document}

\title{A Study of Piecewise Linear-Quadratic Programs
\thanks{The work of the first and forth authors was based on research partially
supported by the U.S.\ National Science Foundation grant IIS--1632971. The work of the second author was supported by the National Natural Science Foundation of
China (Grant No.\ 61731018) and CUHK(SZ) President's Fund (PF.01.000183). The work of the third author was based on
research partially supported by the U.S.\ National Science Foundation grant CMMI-1727757.}}

\titlerunning{A Study of Piecewise Linear-Quadratic Programs}        

\author{Ying Cui  \and Tsung-Hui Chang \and Mingyi Hong \and \\
 Jong-Shi Pang
}

\institute{Ying Cui \at
              The Daniel J.\ Epstein
Department of Industrial and
Systems Engineering, University of Southern California, Los Angeles, CA 90089, U.S.A. \\
 \email{yingcui@usc.edu}
           \and
          Tsung-Hui Chang \at
            School of Science and Engineering, The Chinese University of Hong Kong, Shenzhen, 518172 China \\
              \email{changtsunghui@cuhk.edu.cn}
              \and
              Mingyi Hong \at
             Department of Electrical and Computer Engineering,
University of Minnesota, Minneapolis, Minnesota 55455, U.S.A.\\
\email{mhong@umn.edu}
\and
Jong-Shi Pang \at
 The Daniel J.\ Epstein
Department of Industrial and
Systems Engineering, University of Southern California, Los Angeles, CA 90089, U.S.A. \\
 \email{jongship@usc.edu}
}

\date{June 11, 2018}

\maketitle

\begin{abstract}
Motivated by a growing list of nontraditional statistical estimation problems of the piecewise kind \cite{CuiPangSen18}, this paper provides a survey of known results supplemented with new results for the class of piecewise linear-quadratic programs. These are linearly constrained optimization problems with piecewise linear-quadratic (PLQ) objective functions.  Starting from a study of the representation of such a function in terms of a family of elementary functions consisting of squared affine functions, squared plus-composite-affine functions, and affine functions themselves, we summarize some local properties of a PLQ function in terms of their first and second-order directional derivatives.  We extend some well-known necessary and sufficient second-order conditions for local optimality of a quadratic program to a PLQ program and provide a dozen such equivalent conditions for strong, strict, and isolated local optimality, showing in particular that a PLQ program has the same characterizations for local minimality as a standard quadratic program.   As a consequence of one such condition, we show that the number of strong, strict, or isolated local minima of a PLQ program is finite; this result supplements a recent result about the finite number of directional stationary objective values.  Interestingly, these finiteness results can be uncovered by invoking a very powerful property of subanalytic functions; our proof is fairly elementary, however. We discuss applications of PLQ programs in some modern statistical estimation problems.  
These problems lead to a special class of unconstrained
composite programs involving the non-differentiable $\ell_1$-function, for which we show that the task of verifying the second-order stationary condition can be converted to the problem of checking the copositivity of certain Schur complement on the nonnegative orthant.
\keywords{piecewise linear-quadratic programming \and directional stationarity \and second-order local optimality theory   \and second-order directional, semi- and sub-derivatives \and statistical estimation problems \and matrix copositivity}
\subclass{ 	90C20 \and 	90C26   \and 	68Q25 }
\end{abstract}


\section{Introduction} \label{sec:introduction}

The subject of quadratic programming is as old as that of linear programming.   The monograph \cite{LTYen05}
provides a valuable reference collecting in one source the fundamental theory of quadratic programming.  A piecewise
linear-quadratic (PLQ) function
is a continuous function whose domain is the union of finitely many polyhedral sets on each of which the function is
quadratic.  A piecewise linear-quadratic program is an optimization problem with a PLQ objective and linear constraints.
It appears
that the Ph.D.\ thesis \cite{Sun86} is the first systematic study of a PLQ program; this is followed by the published paper \cite{Sun92}
which studies the class of convex PLQ functions and establishes many structural properties of such piecewise functions.  Like
a quadratic function/program
that provides a bridge between a linear function/program to a nonlinear one, a PLQ function/program provides an important gateway
to a general piecewise smooth function/program from a piecewise affine function/program.  Formal definitions of all these piecewise functions
will be reviewed in Section~\ref{sec:preliminaries}.  A wealth of basic properties of PLQ functions/programs has been obtained in
the treatise \cite{RockafellarWets09}, some of which are succinctly summarized in the most recent
article \cite{BurkeEngle18} and employed in the convergence analysis of Newton and quasi-Newton methods for convex composite programs.
In addition, there are scattered studies of PLQ functions/programs such as the recent one \cite{CuiPang18}
which shows among other things that the set of directional stationary values of the objective function of a PLQ program is finite, in spite
of the possible continuum of local minima of such a problem.   In spite of the abundance of results in the existing studies, there are
some open questions about a PLQ function/program that deserve to be answered.

Motivated by a growing list of applications in several areas, this paper puts together in one place some modern realizations of
PLQ functions/programs, surveys known results to date about these programs, and supplements the old results by new ones.  The
new results aim at addressing some natural questions arising from existing results for quadratic programs and piecewise affine functions.
In the process,
we also clarify some second-order properties of piecewise smooth functions.  We provide here the background
for the new results that we will detail in Section~\ref{sec:preliminaries}.
\\[0.1in]
$\bullet $ On one hand, it is an elementary linear-algebraic fact that a quadratic function is equal to the difference of two sums of squares of affine functions
plus a separate affine function, by the eigen-decomposition of the quadratic form.   On the other hand, it is known from \cite{BartelsKuntzScholtes94,Scholtes02}
that a piecewise affine function admits a max-min representation in terms of affine functions.  This representation is of an algebraic flavor and
is different from the structural properties
of a PLQ function as summarized in \cite[Lemma~2.50]{RockafellarWets09}; see also \cite[Theorem~6.1]{BurkeEngle18}.  Prior results for convex PLQ functions
can be found in \cite{Sun86,Sun92} as mentioned above.  In spite of these known results, there is an absence of an algebraic representation of a PLQ function
that unifies those of a quadratic function and a piecewise affine function.  Part of the contributions of this paper is to provide one such algebraic representation
for a PLQ function, refining the proof of \cite[Proposition~11]{NouiehedPangRazaviyayn17} that provides a difference-of-convex representation of piecewise
functions with a convex domain equal to the union of finitely many closed convex pieces on each of which the function gradients are Lipschitz continuous.
\\[0.1in]
$\bullet $ The study of optimality conditions for constrained optimization problems dates back five decades to the
beginning years of nonlinear programming \cite{FiaccoMcCormick90} under twice continuous differentiability
of the defining functions.  In the early 1980's, such conditions are extended to directionally differentiable
problems using one-sided directional derivatives \cite{BenTalZowe85,BenTalZowe82a,BenTalZowe82b,Chaney87}.  The study of
optimality conditions continues to the modern era of variational analysis \cite{RockafellarWets09}
and generalized differentiation \cite{Mordukhovich06} where the treatment is based on some robust concepts
of first-order {\sl subgradients} and second-order {\sl subderivatives}.  In particular, results from
variational analysis \cite[Theorems 10.1 and 13.24]{RockafellarWets09} establish the necessity of the second-order conditions
for local optimality and the necessity and sufficiency of the strengthened second-order conditions for strong local optimality
for general nonsmooth functions in terms of such subgradients and subderivatives.
Since the early days of quadratic programming \cite{Contesse80,Majthay71}, it is known that the second-order necessary
conditions are indeed sufficient for local optimality and the second-order sufficient conditions are necessary for strong local
optimality.  These results are extended in \cite{Borwein82} to convex constrained
quadratic programs.  Another contribution of this paper is to extend these results in classical quadratic programming to the
class of linearly constrained PLQ programs, thus closing the gap of the local minimality characterizations for this
class of nonsmooth optimization problems.
\\[0.1in]
$\bullet $ Due to the piecewise structure of a PLQ function, it is natural to establish the local minimality of a PLQ in terms of its ``pieces''
which are standard quadratic programs.  In the case of strong local minimality, we provide, via the theory of isolated solutions of affine variational
inequalities \cite[Section~3.3]{FacchineiPang03}, a dozen necessary and sufficient conditions among which are the equivalence of strong, strict, and isolated
local minima \cite{Robinson82} and a matrix-theoretic characterization pertaining to the pieces.  Interestingly, the latter characterization enables us
to show that the number of such minima is finite.  This finiteness result complements similar
results for the objective values of directional stationary solutions; see \cite{CuiPang18}.  As it turns out, these finiteness results for quadratic problems
can be derived by invoking (through additional arguments) a very powerful property of subanalytic sets \cite{BolteDaniilidisLewis06}
whose proof requires advanced mathematical concepts
and abstract analysis.  In contrast, our proof in Proposition~\ref{pr:finite strong} makes use of simple arguments and highlights one consequence of the
necessity of the second-order sufficient conditions for such minima.   The connection between the abstract result in \cite{BolteDaniilidisLewis06}
and our proof also sheds light on the technical difficulty in extending these results to general piecewise quadratic programs whose pieces can be quite
arbitrary.

In addition to these theoretical contributions, we present a host of modern statistical estimation problems that can be formulated as PLQ optimization problems,
and discuss a class of unconstrained
composite programs involving the non-differentiable absolute-value function.  For this special problem, we show that the task of verifying the second-order
stationary condition can be converted to the problem of checking the copositivity of certain Schur complement on the nonnegative orthant.


\section{Preliminaries and Background Results} \label{sec:preliminaries}

Divided into five subsections, this section collects the concepts and background results about
directional derivatives and their role in the optimality conditions of nonsmooth functions as well
as the second-order optimality theory for quadratic programs.  These are summarized here as
a review and also for ease of later reference.  Subsection~\ref{subsec:derivatives} introduces the first- and second-directional
derivatives and related definitions.  Subsection~\ref{subsec:piecewise}
reviews piecewise functions and their local properties and state the max-min representation of a piecewise affine function.
Subsection~\ref{subsec:SC1 functions} discusses the class of semismoothly differentiable (SC$^{\, 1}$) functions which contain the piecewise
smooth functions.  Subsection~\ref{subsec:locmin and stationarity} defines
various local minimizers and first- and second-order stationary points in terms of certain first and second-order necessary
and sufficient conditions.  We also connect these conditions to an abstract result for a general nonsmooth problem.
Subsection~\ref{subsec:QP} summarizes the optimality results for standard quadratic programs.

\subsection{Directional derivatives}  \label{subsec:derivatives}

The following definitions of directional derivatives can all be found in \cite{RockafellarWets09}.
Let $f : \Omega \to \mathbb{R}$ be a given function defined on the open set $\Omega \subseteq \mathbb{R}^n$.
The  (first-order) subderivative ${\rm d}f(x)(v)$ and
one-sided directional derivative $f^{\, \prime}(x;v)$ at a point $x\in \Omega$ along the direction $v\in \mathbb{R}^n$ are
defined by, respectively.
\[
{\rm d}f(x)(v) \, \triangleq \, \displaystyle
\operatornamewithlimits{liminf}_{\substack{v'\to v \\ \tau \downarrow 0}}
\, \displaystyle
\frac{f(x + \tau \, v') - f(x)}{\tau}
\epc \mbox{and}\epc
f^{\, \prime}(x;v) \, \triangleq \, \displaystyle{
\lim_{\tau \downarrow 0}
} \, \displaystyle{
\frac{f(x + \tau \, v) - f(x)}{\tau}
}.
\]
The function $f$ is {\sl directionally differentiable} at $x$ if $f^{\, \prime}(x;v)$ exists for all $v \in \mathbb{R}^n$;
$f$ is {\sl semidifferentiable} at $x$  \cite[Definition~7.20]{RockafellarWets09}
if the ``liminf'' giving ${\rm d}f(x)$ coincides with the ``limsup''; i.e., if the limit
\begin{equation} \label{eq:dd limit}
\displaystyle
\operatornamewithlimits{\lim}_{\substack{v^{\prime} \to v \\ \tau \downarrow 0}}
\, \displaystyle{
\frac{f(x + \tau \, v^{\prime}) - f(x)}{\tau}
}
\end{equation}
exists for all $v \in \mathbb{R}^n$; in this case, we have
\[
{\rm d}f(x)(v) \, = \, \displaystyle
\operatornamewithlimits{\lim}_{\substack{v'\to v \\ \tau \downarrow 0}}
\, \displaystyle{
\frac{f(x + \tau \, v') - f(x)}{\tau}
} \, = \, f^{\, \prime}(x;v), \epc \forall \, v \, \in \, \mathbb{R}^n.
\]
The function $f$ is {\sl B(ouligand) differentiable} at $x$ if it is directionally differentiable at $x$ and locally Lipschitz continuous
near $x$; the latter means that $f$ is Lipschitz continuous in a neighborhood of $x$.
It is easy to see that if $f$ is B-differentiable at $x$,
then the limit (\ref{eq:dd limit}) exists and equals $f^{\, \prime}(x;v)$ for all $v$; moreover,
in this case, the directional derivative $f^{\, \prime}(x;\bullet)$ is Lipschitz continuous on
$\mathbb{R}^n$; see \cite{Shapiro90}.  Thus if $f$ locally Lipschitz continuous near $x$, then semidifferentiability at $x$ is
equivalent to directional differentiability at $x$.

Extending the first-order directional derivative concepts, we define
the {\sl second-order directional derivative} of $f$ at a point $x \in \Omega$ along
the direction $v \in \mathbb{R}^n$ as
\begin{equation} \label{eq:directional second-order dd}
\thalf \, f^{(2)}(x;v) \, \triangleq \, \displaystyle{
\lim_{\tau \downarrow 0}
} \, \displaystyle{
\frac{f(x + \tau \, v) - f(x) - \tau \, f^{\, \prime}(x;v)}{\tau^2},
}
\end{equation}
if the limit exists; the second-order subderivative \cite[Definition~13.3]{RockafellarWets09} at $x$ for $v, w \in \mathbb{R}^n$
is defined as
\[
\thalf \, {\rm d}^2 f(x \,|\, v)(w)\,\triangleq \, \operatornamewithlimits{liminf}_{\substack{w'\to w \\ \tau \downarrow 0}}
\frac{f(x + \tau \, w^{\, \prime}) - f(x) - \tau \, v^Tw^{\, \prime}}{\tau^2}.
\]
Clearly, $f^{(2)}(x;\bullet)$ and ${\rm d}^2 f(x \,|\, v)(\bullet)$ are both positively homogeneous functions of degree 2.
Unlike the directional derivative $f^{\, \prime}(x;\bullet)$
which is a Lipschitz function when $f$ is locally Lipschitz at $x$, the second-order directional derivative $f^{(2)}(x;\bullet)$
is not necessarily continuous; see the PQ function (\ref{eq:PQ example}).
Based on the first-order subderivative ${\rm d}f(x)(\bullet)$, the following second-order subderivative (without mentioning $v$) can be defined:
\[
\thalf \, {\rm d}^2 f(x)(w)\,\triangleq \, \displaystyle{
\liminf_{\substack{w^{\, \prime} \to w \\ \tau \downarrow 0}}
} \, \displaystyle{
\frac{f(x + \tau \, w^{\, \prime}) - f(x) - \tau \, {\rm d}f(x)(w^{\, \prime})}{\tau^2}
} \, .
\]
We say that $f$ is {\sl twice directionally differentiable} at $x$ if it is directionally differentiable at $x$ and the limit $f^{(2)}(x;v)$ exists
for all $v \in \mathbb{R}^n$.
According to \cite[Definition~13.6]{RockafellarWets09}, $f$ is said to be {\sl twice semidifferentiable} at $x$ if it is semidifferentiable at $x$
and the limit
\begin{equation} \label{eq:2-order semidiff limit}
\displaystyle{
\lim_{\substack{w^{\, \prime} \to w \\ \tau \downarrow 0}}
} \, \displaystyle{
\frac{f(x + \tau \, w^{\, \prime}) - f(x) - \tau \, {\rm d}f(x)(w^{\, \prime})}{\tau^2}
}
\end{equation}
exists for all $w \in \mathbb{R}^n$.  If $f$ is 
twice semidifferentiable at $x$, then for all $w \in \mathbb{R}^n$,
\[
\thalf {\rm d}^2 f(x)(w) \, = \, \displaystyle{
\lim_{\substack{x^{\, \prime} \to x \\ \tau \downarrow 0}}
} \, \left\{ \, \displaystyle{
\frac{f(x^{\, \prime}) - f(x) - f^{\, \prime}(x;x^{\, \prime} - x)}{\tau^2}
} \, : \, \displaystyle{
\frac{x^{\, \prime} - x}{\tau}
} \, \to w \, \right\} \, = \, \thalf \, f^{(2)}(x;w) .
\]
Moreover, in this case, ${\rm d}^2 f(x)(\bullet)$, and thus $f^{(2)}(x;\bullet)$, is continuous, by \cite[Exercise~13.7]{RockafellarWets09}.

\subsection{Piecewise functions} \label{subsec:piecewise}

We recall that a function
$f$ is PC$^{\, k}$ on an open subset $\Omega$ of $\mathbb{R}^n$ for a positive integer $k$ if it is continuous and there exist finitely
many C$^{\, k}$ (for $k$-times continuously differentiable) functions $\{ \, f_i \, \}_{i=1}^I$
such that $f(x) \in \{ \, f_i(x) \, \}_{i=1}^I$ for all $x \in \Omega$.  For a given $x \in \Omega$, let ${\cal A}(x) \subseteq \{ 1, \cdots, I \}$
be the index set consisting of indices $i$ such that $f(x) = f_i(x)$.  For each $i = 1, \cdots, I$, the pair $(f_i,\Omega^{\, i})$, where
$\Omega^{\, i} \triangleq \left\{ x \in \Omega \mid f(x) = f_i(x) \right\}$ is a called a {\sl piece} of $f$.  Occasionally, we will also call each
function $f_i$ and set $\Omega^{\, i}$ separately a piece of $f$.

Of particular interest in this paper are several classes of piecewise
functions.  We say that a continuous function $f : {\cal D} \to \mathbb{R}$
defined on a set ${\cal D} \subseteq \mathbb{R}^n$ is {\sl piecewise quadratic} (PQ) if there exist finitely many quadratic functions $\{ \, q_i \, \}_{i=1}^I$
such that $f(x) \in \{ \, q_i(x) \, \}_{i=1}^I$ for all $x \in {\cal D}$.  The continuous function $f : {\cal D} \to \mathbb{R}$
is {\sl piecewise linear-quadratic} (PLQ) \cite[Chapters~10.E and 11.D]{RockafellarWets09} if there exist finitely many quadratic functions $\{ q_i \}_{i=1}^I$
and the same number of polyhedra $\{ P^{\, i} \}_{i=1}^I$ whose union is ${\cal D}$
such that $f(x) = q_i(x)$ for all $x \in P^{\, i}$; thus ${\cal D}$ is a closed set.
In the terminology of the cited reference, ${\cal D}$ is called the {\sl domain} of the PLQ function $f$ and is denoted $\mbox{dom }f$.
We call a set $S \subseteq \mathbb{R}^n$ {\sl piecewise polyhedral} if it is the union of finitely many polyhedra each of which is called
a {\sl (polyhedral) piece} of $S$.   Thus the domain of a PLQ function is piecewise polyhedral.
Piecewise quadratic functions need not be piecewise linear-quadratic because there is no requirement for the existence of a family of polyhedral decomposition
of the domain as required by a PLQ function.  A piecewise affine (PA) function is a PLQ function such that the quadratic element functions $q_i$
are all affine functions.

It is known that PC$^{\, 1}$, and thus PQ, functions are B-differentiable;
see e.g.\ Lemma~4.6.1 in \cite{FacchineiPang03}.  Moreover, the directional derivative $f^{\, \prime}(x;d)$ is equal to $\nabla f_i(x)^Td$ for every
index $i \in {\cal A}^{\, \prime}(x;d)$, where ${\cal A}^{\, \prime}(x;d)$, called the {\sl directionally active set} at $x$ in the direction $d$,
consists of those indices $i^{\, \prime}$ for which there
exists a sequence of positive scalars $\{ \tau_{\nu} \}$ converging to zero such that $f_{i^{\, \prime}}(x + \tau_{\nu} d) = f(x + \tau_{\nu} d)$
for all $\nu$.  Implicit in this result is the fact that $\nabla f_i(x)^Td = \nabla f_j(x)^Td$ for any two indices $i$ and $j$ in ${\cal A}^{\, \prime}(x;d)$.
Thus the directional derivative $f^{\, \prime}(x;\bullet)$ of a PC$^{\, 1}$ function is a piecewise linear function on $\mathbb{R}^n$.  A generalization of this result
is proved for PC$^{\, 2}$ functions in Proposition~\ref{pr:PC are dd2} that extends the result below for PLQ functions; a remark following the latter proposition
highlights the difference between twice directional differentiability and twice semidifferentiability.
In the following result and subsequently, ${\cal T}(\bar{x};S)$ denotes the tangent cone of a closed set $S$ at a point $\bar{x} \in S$; i.e.,
$v \in {\cal T}(\bar{x};S)$ if and only if there exist a sequence of vectors $\{ x^k \} \subset S$ converging to $\bar{x}$ and a sequence of positive
scalars $\{ \tau_k \} \downarrow 0$ such that $v = \displaystyle{
\lim_{k \to \infty}
} \, \displaystyle{
\frac{x^k - \bar{x}}{\tau_k}
}$.

\begin{proposition}\label{prop:convex plq}\rm \cite[Proposition 13.9]{RockafellarWets09}
Let $f: {\cal D} \subseteq \mathbb{R}^n \to \mathbb{R}$ be a PLQ function with the domain
${\cal D}$ being the union of the polyhedral pieces $\{ P^{\, i} \}_{i=1}^I$; associated with each of such piece $P^{\, i}$
is the quadratic function $q_i$ for $i = 1, \cdots, I$.  At any point
$\bar{x} \in {\rm dom}\,f$, $f^{\, \prime}(\bar{x};\,\bullet\,) = {\rm d}f(\bar{x})$, which is piecewise linear with
${\rm dom}\,{\rm d}f(\bar{x}) = \mathcal{T}(\bar{x}; {\rm dom}\,f)$.  In particular, for $i \in  \mathcal{A}(\bar{x})$ and
$v \in \mathcal{T}(\bar{x}; P^{\,i})$,
\[
f^{\, \prime}(\bar{x};v) \, = \, \nabla q_i(\bar{x})^T v . \]
In addition,
$f^{(2)}(\bar{x};\,\bullet\,) = {\rm d}^2 f(\bar{x})$ is piecewise linear-quadratic given by
\[
f^{(2)}(\bar{x};v) = {\rm d}^2 f(\bar{x})(v) =
\left\{\begin{array}{ll}
v^T \nabla^2 q_i(\bar{x}) v & \;\,\mbox{if $v\in \mathcal{T}(\bar{x}; P^{\,i})$} \\[0.05in]
+\infty & \;\,\mbox{otherwise}.
\end{array}
\right.
\]
Moreover, there exists a neighborhood $\mathcal{N}$ of $\bar{x}$ such that
\[ f(x) = f(\bar{x}) + f^{\, \prime}(\bar{x}; x- \bar{x}) +f^{(2)}(\bar{x};x- \bar{x}),  \epc \forall \; x \in {\rm dom}\,f \, \cap\, \mathcal{N}.
\]
As noted in \cite[Proposition~4.2]{BurkeEngle18}, no convexity on $f$ is needed in the above statements.  \hfill $\Box$.
\end{proposition}

According to \cite{BartelsKuntzScholtes94,Scholtes02}, every PA function with domain $\mathbb{R}^n$ has a max-min representation.  Specifically, if
$f : \mathbb{R}^n \to \mathbb{R}$ is PA, then there exist finitely many affine functions $\left\{ \, \left( f_{ij} \right)_{j=1}^{J_i} \, \right\}_{i=1}^I$
such that
\begin{equation} \label{eq:Scholtes PA}
f(x) \, = \, \displaystyle{
\max_{1 \leq i \leq I}
} \, \displaystyle{
\min_{1 \leq j \leq J_i}
} \, f_{ij}(x), \epc x \, \in \, \mathbb{R}^n.
\end{equation}
From this representation, it is easy to deduce that if $f$ is PA, then
\[
f(x^{\, \prime}) \, = \,  f(x) + f^{\, \prime}(x;x^{\, \prime} -x), \epc \forall \, x^{\, \prime} \mbox{ near } x.
\]
From the representation (\ref{eq:Scholtes PA}), we may deduce that every PA function is a difference-of-convex (dc) function with the following
difference-max-affine representation:
\begin{equation} \label{eq:dc representation PA}
f(x) \, = \, \underbrace{\displaystyle{
\max_{1 \leq i \leq I_1}
} \, \left( \, x^T a^i + \alpha_i \, \right)}_{\mbox{convex in $x$}} -
\underbrace{\displaystyle{
\max_{1 \leq i \leq I_2}
} \, \left( \, x^T b^i + \beta_i \, \right)}_{\mbox{convex in $x$}}
\end{equation}
for some positive integers $I_1$ and $I_2$, $n$-vectors $\{ a^i \}_{i=1}^{I_1}$ and $\{ b^i \}_{i=1}^{I_2}$, and scalars $\{ \alpha_i \}_{i=1}^{I_1}$
and $\{ \beta_i \}_{i=1}^{I_2}$.   In view of the two algebraic representations (\ref{eq:Scholtes PA}) and (\ref{eq:dc representation PA}), it is
natural to ask whether a PLQ function has similar representations using quadratic functions.  This question easily has a negative answer as illustrated
by the squared plus function; i.e., $t_+^2$, where $t_+ \triangleq \max(t,0)$.  Incidentally, the latter representation (\ref{eq:dc representation PA})
is key to the statistical estimation problem using a PA model; see \cite{CuiPangSen18}.  By a result in the recent paper \cite{NouiehedPangRazaviyayn17}, which we
rephrase below, it follows that that every piecewise quadratic function with a convex domain is a dc function.  A function is LC$^{\, 1}$ if it is differentiable
with a Lipschitz gradient.  No convexity of the function $\theta$ is required in the proposition.

\begin{proposition} \label{pr:Maher PQ} \rm \cite[Proposition~11]{NouiehedPangRazaviyayn17}
Let $\theta(x)$ be a continuous function on a convex set
${\cal S} \triangleq \displaystyle{
\bigcup_{i=1}^I
} \, S^{\, i}$ where each $S^{\, i}$ is a closed convex set in $\mathbb{R}^N$.  Suppose there exist LC$^{\, 1}$
functions $\{ \theta_i(x) \}_{i=1}^I$ defined on an open set ${\cal O}$ containing ${\cal S}$
such that $\theta(x) = \theta_i(x)$ for all $x \in S^{\, i}$ and that each difference function
$\theta_{ji}(x) \triangleq \theta_j(x) - \theta_i(x)$ has dc gradients on ${\cal S}$.
It holds that $\theta$ is dc on ${\cal S}$ with the following representation:
\begin{equation} \label{eq:Maher PQ cd representation}
\theta(x) \, = \, \displaystyle{
\min_{1 \leq i \leq I}
} \, \left\{ \, \theta_i(x) + \mbox{dist}_2(x;S^{\, i}) \, \displaystyle{
\max_{1 \leq j \leq I}
} \, \| \nabla \theta_{ji}(x) \|_2 + \displaystyle{
\frac{3 \, L_i}{2}
} \, \left[ \, \mbox{dist}_2(x;S^{\, i}) \, \right]^2 \, \right\} \epc \forall \, x \, \in \, {\cal S},
\end{equation}
where $\mbox{dist}_2(x;S^{\, i}) \triangleq \displaystyle{
\operatornamewithlimits{\mbox{minimum}}_{y \in S_i}
} \, \| \, y - x \, \|_2$ is the Euclidean distance from $x$ to the set $S^{\, i}$ and the constant
$L_i \triangleq \displaystyle{
\max_{1 \leq j \leq I}
} \, L_{ji}$ with each $L_{ji}$ being a Lipschitz constant of $\nabla \theta_{ji}$.  \hfill $\Box$
\end{proposition}

This result is the starting point to derive an algebraic representation of a PLQ function in terms of some elementary functions.

\subsection{Semismoothly differentiable functions} \label{subsec:SC1 functions}

Piecewise C$^{\, k}$ functions are a subclass of the class of semismooth functions formally defined as follows.
A vector function $\Phi : \Omega \to \mathbb{R}^m$ defined on the open set $\Omega \subseteq \mathbb{R}^n$ is
{\sl semismooth} \cite{FacchineiPang03,Kummer88,Mifflin77,QiSun93} at $\bar{w} \in \Omega $ if $\Phi$
is B-differentiable near $\bar{w}$ and
\[
\displaystyle{
\lim_{\stackrel{\bar{w} \neq w \to \bar{w}}{H \in \partial \Phi(w)}}
} \, \displaystyle{
\frac{\Phi^{\prime}(\bar{w};w - \bar{w}) - H \, ( \, w - \bar{w} \, )}{\| \, w - \bar{w} \, \|}
} \, = \, 0,
\]
where $\partial \Phi(w)$ denotes the (generalized) Clarke Jacobian \cite{Clarke83} of $\Phi$ at $w$.
A continuous real-valued function $\psi : \Xi \subseteq \mathbb{R}^m \to \mathbb{R}$ defined on the open set $\Xi$ is
{\sl semismoothly differentiable} (SC$^{\, 1}$) at $\bar{z} \in \Xi$
if it is once differentiable near $\bar{z}$ and its gradient is semismooth at $\bar{z}$.
By \cite[Proposition~7.4.10; expression (7.4.14) more precisely]{FacchineiPang03}, it holds that
if $\psi$ is SC$^{\, 1}$ at $\bar{z}$ with semismooth gradient $\Psi$, then
\begin{equation} \label{eq:SC1 limit}
\displaystyle{
\lim_{d \to 0}
} \, \displaystyle{
\frac{\psi(\bar{z} + d) - \psi(\bar{z}) - \Psi(\bar{z})^Td - \thalf \, d^{\, T} \Psi^{\, \prime}(\bar{z};d)}{\| \, d \, \|^2}
} \, = \, 0.
\end{equation}
The next result shows in particular that a SC$^{\, 1}$ function must be twice semidifferentiable.  This result adds a new local
property of a SC$^{\, 1}$ function.  See Section~\ref{sec:a class of composite} for an application of the result.

\begin{proposition} \label{pr:SC1 is dd2} \rm
Let $f : \Xi \subseteq \mathbb{R}^m \to \mathbb{R}$ be SC$^{\, 1}$ near $\Phi(\bar{w}) \in \Xi$
and $\Phi : \Omega \subseteq \mathbb{R}^n \to \Xi$ be locally Lipschitz and twice semidifferentiable near
$\bar{w} \in \Omega$.
The composite function $\varphi \triangleq f \circ \Phi : \Omega \to \mathbb{R}$
is  twice semidifferentiable at $\bar{w}$; moreover, with $F(y) \triangleq \nabla f(y)$
\begin{equation} \label{eq:dd2 composite}
\varphi^{\, (2)}(\bar{w};v) \, = \,
\Phi^{\, \prime}(\bar{w};v)^T F^{\, \prime}(\Phi(\bar{w});\Phi^{\, \prime}(\bar{w};v)) + F(\Phi(\bar{w}))^T \Phi^{(2)}(\bar{w};v),
\epc \mbox{for all $v \in \mathbb{R}^n$}.
\end{equation}
\end{proposition}
{\bf Proof.}  It suffices to show that
\[
\displaystyle{
\lim_{\substack{v^{\, \prime} \to v \\ \tau \downarrow 0}}
} \, \displaystyle{
\frac{\varphi(\bar{w} + \tau \, v^{\, \prime}) - \varphi(\bar{w}) - \tau \, \varphi^{\, \prime}(\bar{w};v^{\, \prime})}{\tau^2}
} \, = \, \Phi^{\, \prime}(\bar{w};v)^T F^{\, \prime}(\Phi(\bar{w});\Phi^{\, \prime}(\bar{w};v)) + F(\Phi(\bar{w}))^T \Phi^{(2)}(\bar{w};v).
\]
Since $\varphi^{\, \prime}(\bar{w};v^{\, \prime}) = F(\Phi(\bar{w}))^T\Phi^{\, \prime}(\bar{w};v^{\, \prime})$, writing
$d\Phi \triangleq \Phi(\bar{w} + \tau \, v^{\, \prime}) - \Phi(\bar{w})$, we have
\[ \begin{array}{l}
\displaystyle{
\frac{\varphi(\bar{w} + \tau \, v^{\, \prime}) - \varphi(\bar{w}) - \tau \, \varphi^{\, \prime}(\bar{w};v^{\, \prime})}{\tau^2}
} \\ [0.2in]
\epc = \, \displaystyle{
\frac{f(\Phi(\bar{w} + \tau \, v^{\, \prime})) - f(\Phi(\bar{w})) - F(\Phi(\bar{w}))^Td\Phi - \thalf \, d\Phi^T F^{\, \prime}(\bar{w};d\Phi)}{
\| \, \Phi(\bar{w} + \tau \, v^{\, \prime}) - \Phi(\bar{w}) \, \|^2} } \, \displaystyle{
\frac{\| \, \Phi(\bar{w} + \tau \, v^{\, \prime}) - \Phi(\bar{w}) \, \|^2}{\tau^2}
} \\ [0.2in]
\hspace{0.5in} + \, \thalf \, \displaystyle{
\frac{\left[ \, \Phi(\bar{w} + \tau \, v^{\, \prime}) - \Phi(\bar{w}) \, \right]^TF^{\, \prime}(\bar{w};\Phi(\bar{w} + \tau \, v^{\, \prime})
- \Phi(\bar{w}))}{\tau^2}
} \\ [0.2in]
\hspace{0.5in} + \, \displaystyle{
\frac{F(\Phi(\bar{w}))^T \left[ \, \Phi(\bar{w} + \tau \, v^{\, \prime}) - \Phi(\bar{w}) - \tau \, \Phi^{\, \prime}(\bar{w};v^{\, \prime}) \, \right]}{\tau^2}
} \, .
\end{array}
\]
Since
\[
\displaystyle{
\lim_{\substack{v^{\, \prime} \to v \\ \tau \downarrow 0}}
} \, \displaystyle{
\frac{\Phi(\bar{w} + \tau \, v^{\, \prime}) - \Phi(\bar{w})}{\tau}
} \, = \, \Phi^{\, \prime}(\bar{w};v) \mbox{ and }
\displaystyle{
\lim_{\substack{v^{\, \prime} \to v \\ \tau \downarrow 0}}
} \, \displaystyle{
\frac{\Phi(\bar{w} + \tau \, v^{\, \prime}) - \Phi(\bar{w}) - \tau \, \Phi^{\, \prime}(\bar{w};v^{\, \prime})}{\thalf \, \tau^2}
} \, = \, \Phi^{(2)}(\bar{w};v),
\]
combining these limits with (\ref{eq:SC1 limit}) applied to $f$ with gradient $F$ at $\Phi(\bar{w})$, we easily obtain the desired
formula (\ref{eq:dd2 composite}) for $\varphi^{\, (2)}(\bar{w};v)$.  \hfill $\Box$

\subsection{Local minimizers and stationarity} \label{subsec:locmin and stationarity}

Consider the optimization problem:
\begin{equation} \label{eq:basic opt problem}
\displaystyle{
\operatornamewithlimits{\mbox{minimize}}_{x \in X}
} \ f(x),
\end{equation}
where $X$ is a polyhedral set (unless otherwise specified) in $\mathbb{R}^n$ and $f$ is a locally Lipschitz continuous function defined
on an open set containing $X$.  We say that $\bar{x} \in X$ is a

\noindent
$\bullet $ {\bf local minimizer} of $f$ on $X$ if there exists an (open) neighborhood ${\cal N}$ of $\bar{x}$ such that
$f(x) \geq f(\bar{x})$ for all $x \in X \cap {\cal N}$;
\\[0.05in]
$\bullet $ {\bf strict local minimizer} of $f$ on $X$ if there exists an (open) neighborhood ${\cal N}$ of $\bar{x}$ such that
$f(x) > f(\bar{x})$ for all $x \in X \cap {\cal N}$ and $x \neq \bar{x}$;
\\[0.05in]
$\bullet $ {\bf isolated local minimizer} of $f$ on $X$ if there exists an (open) neighborhood ${\cal N}$ of $\bar{x}$ such that
$\bar{x}$ is the only local minimizer in ${\cal N}$ of $f$ constrained by $X$;
\\[0.05in]
$\bullet $ {\bf strong local minimizer} of $f$ on $X$ if there exist a scalar $c > 0$ and an (open) neighborhood ${\cal N}$ of $\bar{x}$ such that
$f(x) \geq f(\bar{x}) + c \| \, x - \bar{x} \, \|^2$ for all $x \in X \cap {\cal N}$.

Clearly every strong local minimizer must be strict; so is every isolated local minimizer.  It is known that the converse of these
statements are not valid for a general nonlinear program.
Stated for a proper extended-valued function, i.e., $f \not\equiv \infty$,
the following theorem provides a general result for the local optimality based on the first and second-order subderivatives.

\begin{theorem}\label{thm: RW: opt}\rm  \cite[Theorem 10.1 \& 13.24]{RockafellarWets09}
Let $f:\mathbb{R}^n\to (-\infty, +\infty]$ be a proper extended-valued function.  The following two statements (a) and (b) hold for the program:
\[
\operatornamewithlimits{minimize}_{x\in \mathbb{R}^n}\; f(x),
\]
(a) If $\bar{x}$ is a local minimum, then ${\rm d}f(\bar{x})(v)\geq 0$ and ${\rm d}^2 f(\bar{x}\,|\,0)(v) \geq 0$ for any $v\in \mathbb{R}^n$.\\[0.05in]
(b)  $\bar{x}$ is a strong local minimum solution if and only if ${\rm d}f(\bar{x})(v)\geq 0$ and ${\rm d}^2 f (x \,|\, 0)(v) > 0$ for all $v\neq  0$.
\hfill $\Box$
\end{theorem}

To apply the above theorem to the problem (\ref{eq:basic opt problem}), one needs to employ the indicator function of the constraint set $X$,
defined as $\delta_X(x) = \left\{ \begin{array}{ll}
0 & \mbox{if $x \in X$} \\
\infty & \mbox{otherwise},
\end{array} \right.$ and
form the extended-valued function $\wh{f}(x) = f(x) + \delta_X(x)$.  With the goal of exposing the constraint
set $X$ in the optimality conditions and avoiding the definition of second-order tangent sets \cite[Section~11.C]{RockafellarWets09}
\cite[Section~3.2.1]{BonnansShapiro00}, which incidentally may not be needed because of the polyhedrality of $X$,
we bypass this extended-valued maneuver and present the following variant
of Theorem~\ref{thm: RW: opt}.  We offer a detailed proof of the implication (b3) $\Rightarrow$ (b1) in the proposition because
we cannot identify a result in the literature that we cite directly.

\begin{proposition} \label{pr:second-order optimality} \rm
Let $f : \Omega \to \mathbb{R}$ be locally Lipschitz continuous near a given $\bar{x} \in X$ and twice semidifferentiable at $\bar{x}$,
where $X$ is a polyhedron contained in the open set $\Omega$.  Consider two sets of statements for the program (\ref{eq:basic opt problem})
at $\bar{x} \in X$.\\[0.05in]
(a1) $\bar{x}$ is a local minimizer;\\[0.05in]
(a2) ${\rm d}f(\bar{x})(v) =  f^{\, \prime}(\bar{x};v) \geq 0$ for all $v \in {\cal T}(\bar{x};X)$, and ${\rm d}^2 f(\bar{x})(v) \geq 0$
for all $v \in {\cal T}(\bar{x};X)$ such that ${\rm d}f(\bar{x})(v) = 0$;\\[0.05in]
(a3) ${\rm d}f(\bar{x})(x - \bar{x}) = f^{\, \prime}(\bar{x};x - \bar{x}) \geq 0$ for all $x \in X$, and ${\rm d}^2 f(\bar{x})(x - \bar{x}) \geq 0$
for all $x \in X \setminus \{ \bar{x} \}$ such that ${\rm d}f(\bar{x})(x - \bar{x}) = 0$;\\[0.05in]
(b1) $\bar{x}$ is a strong local minimizer;\\[0.05in]
(b2) ${\rm d}f(\bar{x})(v) \geq 0$ for all $v \in {\cal T}(\bar{x};X)$, and ${\rm d}^2 f(\bar{x})(v) > 0$
for all nonzero $v \in {\cal T}(\bar{x};X)$ such that ${\rm d}f(\bar{x})(v) = 0$;\\[0.05in]
(b3) ${\rm d}f(\bar{x})(x - \bar{x}) \geq 0$ for all $x \in X$, and ${\rm d}^2 f(\bar{x})(x - \bar{x}) > 0$
for all $x \in X \setminus \{ \bar{x} \}$ such that ${\rm d}f(\bar{x})(x - \bar{x}) = 0$.\\[0.05in]
It holds that (b1) $\Leftrightarrow$ (b2) $\Leftrightarrow$ (b3)  $\Rightarrow$ (a1) $\Rightarrow$ (a2) $\Leftrightarrow$ (a3).
\end{proposition}
{\bf Proof.}  (b1) $\Rightarrow$ (b2).  By the polyhedrality of $X$, it follows that for every $v \in {\cal T}(\bar{x};X)$, $\bar{x} + \tau v \in X$
for all $\tau > 0$ sufficiently small.  Hence the claimed implication is immediate from the equality ${\rm d}^2 f(\bar{x})(v) = f^{(2)}(\bar{x};v)$.

(b2) $\Leftrightarrow$ (b3).  This is easy because $X$ is polyhedral.

(b3) $\Rightarrow$ (b1).  This is nontrivial yet not difficult part of the result.
Assume by way of contradiction
that $\bar{x} \in X$ is not a strong local minimizer.  It then follows that there exists
a sequence $\{ x^k \} \subset X$ converging to $\bar{x}$ such that
\begin{equation} \label{eq:reverse strong locmin}
f(\bar{x}) \, > \, f(x^k) - \displaystyle{
\frac{1}{k}
} \, \| \, \bar{x} - x^k \, \|^2, \epc \forall \, k.
\end{equation}
This implies in particular that $x^k \neq \bar{x}$ for all $k$.  With no loss of generality, we may assume that the normalized sequence
$\left\{ \, \displaystyle{
\frac{x^k - \bar{x}}{\| \, x^k - \bar{x} \, \|}
} \, \right\}$ converges to a limit $v$ which must be nonzero.  Thus, by the continuity of $f^{\, \prime}(\bar{x};\bullet)$,
it follows that
$f^{\, \prime}(\bar{x};v) \geq 0$.  By the local Lipschitz continuity of $f$, we have
\[ \displaystyle{
\lim_{k \to \infty}
} \, \displaystyle{
\frac{f(x^k) - f(\bar{x})}{\| \, x^k - \bar{x} \|}
} \, = \, f^{\, \prime}(\bar{x};v).
\]
Hence (\ref{eq:reverse strong locmin}) yields $f^{\, \prime}(\bar{x};v) \leq 0$.  Thus,
$f^{\, \prime}(\bar{x};v) = 0$.   It follows that ${\rm d}^2 f(\bar{x})(v) > 0$ because $\bar{x} + \tau v \in X$
for all $\tau > 0$ sufficiently small by the polyehdrality of $X$.  Since
\[
\thalf \, {\rm d}^2 f(\bar{x})(v) \, = \, \displaystyle{
\lim_{k \to \infty}
} \, \displaystyle{
\frac{f(x^k) - f(\bar{x}) - f^{\, \prime}(\bar{x};x^k - \bar{x})}{\| \, x^k - \bar{x} \, \|^2},
}
\]
it follows that for some constant $c > 0$,
\[
f(x^k) \, \geq \, f(\bar{x}) + f^{\, \prime}(\bar{x};x^k - \bar{x}) + c \, \| \, x^k - \bar{x} \, \|^2
\, \geq \, f(\bar{x}) + c \, \| \, x^k - \bar{x} \, \|^2.
\]
But this contradicts (\ref{eq:reverse strong locmin}).

The remaining implications (b1) $\Rightarrow$ (a1) $\Rightarrow$ (a2) $\Leftrightarrow$ (a3) are all fairly easy.
\hfill $\Box$
%
%
%

\begin{remark} \rm
Notice that the implication (a3) $\Rightarrow$ (a1) is left out in Proposition~\ref{pr:second-order optimality}.
Inspired by classic results for standard quadratic programs, it
is natural to ask whether for the program (\ref{eq:basic opt problem}) such a reverse implication will be valid if the objective function $f$ is PLQ.
Completing the equivalence of (a1), (a2), and (a3) for a PLQ program is a contribution of this paper.  \hfill $\Box$
\end{remark}

Based on the above result, we define two types of second-order stationary solutions for the problem (\ref{eq:basic opt problem}) with
a twice semidifferentiable function $f$ and a polyhedral $X$ using ${\rm d}^2 f(\bar{x}) = f^{(2)}(\bar{x};\bullet)$.
Specifically, we say that $\bar{x} \in X$
\\[0.05in]
$\bullet $ is a {\bf (directional) stationary point}, or equivalently, satisfies the {\bf (first-order directional) stationarity condition}
if ${\rm d}f(\bar{x})(v) = f^{\, \prime}(\bar{x};v) \geq 0$ for all $v \in {\cal T}(\bar{x};X)$, or equivalently,
${\rm d}f(\bar{x})(x - \bar{x})\geq 0$ for all $x \in X$;
\\[0.05in]
$\bullet $ is an {\bf isolated (or locally unique) stationary point} if there exists an (open) neighborhood
${\cal N}$ of $\bar{x}$ such that $\bar{x}$ is the only stationary point in ${\cal N}$;
\\[0.05in]
$\bullet $ satisfies the {\bf second-order necessary condition} if it is a stationary point and
${\rm d}^2 f(\bar{x})(v) \geq 0$ for all $v \in {\cal T}(\bar{x};X)$ such that ${\rm d}f(\bar{x})(v) = 0$;
\\[0.05in]
$\bullet $ satisfies the {\bf second-order sufficient condition} if it is a stationary point and
${\rm d}^2 f(\bar{x})(v) > 0$ for all nonzero $v \in {\cal T}(\bar{x};X)$ such that ${\rm d}f(\bar{x})(v) = 0$.

Like the first-order stationarity conditions, we also call the second-order necessary and sufficient conditions
second-order stationarity conditions.
Clearly, a local minimizer that is an isolated stationary point must be an isolated local minimizer.  If $\bar{x}$ is a
(directional) stationary point of (\ref{eq:basic opt problem}), we call $f(\bar{x})$ is {\bf (directional) stationary value}
of this problem.

\subsection{Quadratic programs} \label{subsec:QP}

Consider the standard quadratic program:
\begin{equation} \label{eq:QP}
\displaystyle{
\operatornamewithlimits{\mbox{minimize}}_{x \in P}
} \ q(x),
\end{equation}
where $q(x) = \thalf x^TQx + c^Tx + \alpha$ is a quadratic function with the matrix $Q \in \mathbb{R}^{n \times n}$ being symmetric and the pair
$(c,\alpha) \in \mathbb{R}^{n+1}$,
and $P \triangleq \left\{ x \in \mathbb{R}^n \mid Ax \geq b \right\}$ for some matrix $A \in \mathbb{R}^{m \times n}$ and $m$-vector $b$ is a polyhedral set.
We recall that the critical cone, denoted ${\cal C}(\bar{x};q;P)$, of this problem at a given $\bar{x} \in P$ is by definition the
polyhedral cone: ${\cal C}(\bar{x};q;P) \triangleq {\cal T}(\bar{x};P) \, \cap \, \nabla q(\bar{x})^{\perp}$, where $a^{\perp}$ denotes the orthogonal complement of
the vector $a$ consisting of all vectors $v$ perpendicular to $a$.  There is an equivalent definition of the critical cone when the base vector $\bar{x}$ is a
stationary solution of (\ref{eq:QP}) in terms of the constraint multipliers \cite[Section~3.3.1]{FacchineiPang03}.  Specifically, for such a stationary solution
$\bar{x}$, let $\Lambda(\bar{x})$ denote the set of multipliers $\lambda \in \mathbb{R}^m_+$ such that the following Karush-Kuhn-Tucker (KKT) conditions hold:
\[
\begin{array}{lll}
0 & = & \nabla q(\bar{x}) - A^T \lambda \\ [5pt]
0 & \leq & \lambda \, \perp \, A\bar{x} - b \, \geq \, 0,
\end{array} \]
where the $\perp$ notation here denote the complementary slackness property between the constraint multiplier $\lambda$ and the
(nonnegative) slack variable $s \geq A\bar{x} - b$.  Let $\mbox{supp}(\lambda)$ denote the {\sl support} of the vector $\lambda$; i.e., $\mbox{supp}(\lambda)$
consists of all the indices $i \in \{ 1, \cdots, m \}$ such that $\lambda_i > 0$.  We then have
\[
{\cal C}(\bar{x};q;P) \, = \, \left\{ \, v \, \in \, {\cal T}(\bar{x};P) \, \mid \, \exists \lambda \, \in \, \Lambda(\bar{x})
\mbox{ such that } A_{\, i \bullet}v \, = \, 0 \mbox{ for all } i \, \in \, \mbox{supp}(\lambda) \, \right\},
\]
where $A_{\, i \bullet}$ denotes the $i$th row of $A$.  The following result about local minimizers is classical in the theory of quadratic programs.

\begin{proposition} \label{pr:QP local} \rm
A feasible vector $\bar{x} \in P$ of the quadratic program (\ref{eq:QP})is a local minimizer if and only if it satisfies the second-order necessary
condition; this is equivalent to $\bar{x}$ being a stationary point and $Q$ being copositive on ${\cal C}(\bar{x};q;P)$; i.e., $v^TQv \geq 0$
for all $v \in  {\cal C}(\bar{x};q;P)$.  \hfill $\Box$
\end{proposition}

Theorems 3 in \cite{Borwein82} offers an extension of the above proposition to a non-polyhedral feasible set $P$; in the sufficiency part,
it requires the polyhedrality of the tangent cone ${\cal T}(\bar{x};P)$.

The next proposition about strong local minimizers collects various known results from the literature and put them in one place for clarity
and also for ease of later reference; there are a few parts that are not particularly well known but are needed to complete the proof of all
the equivalences.  The $\perp$ notation in part (a5) denotes the orthogonality of the two vectors $v$ and $Qv$;
${\cal C}(\bar{x};q;P)^*$ denotes the dual of the critical cone.

\begin{proposition} \label{pr:QP strong local} \rm
The following statements are equivalent for a feasible vector $\bar{x} \in P$ of the quadratic program (\ref{eq:QP}).\\[0.05in]
(a1) $\bar{x}$ is a stationary point and $Q$ is strictly copositive on ${\cal C}(\bar{x};q;P)$;\\[0.05in]
(a2) $\bar{x}$ is a stationary point, $Q$ is copositive on ${\cal C}(\bar{x};q;P)$, and the implication below holds:
\begin{equation} \label{eq:R0}
\left[ \, {\cal C}(\bar{x};q;P) \, \ni \, v \, \perp \, Qv \, \in \, {\cal C}(\bar{x};q;P)^* \, \right] \ \Rightarrow \ v \, = \, 0;
\end{equation}
(a3) $\bar{x}$ is both a local minimizer and an isolated stationary point;\\[0.05in]
(a4) $\bar{x}$ is an isolated local minimizer;\\[0.05in]
(a5) $\bar{x}$ is a strict local minimizer;\\[0.05in]
(a6) $\bar{x}$ is a strong local minimizer.
\end{proposition}
{\bf Proof.}  The proof follows the implications below which are either easy or known;
\[
\mbox{(a1)} \, \Rightarrow \, \mbox{(a2)} \, \Rightarrow \, \mbox{(a3)} \, \Rightarrow \, \mbox{(a4)} \, \Rightarrow \, \mbox{(a5)}
\, \Leftrightarrow \, \mbox{(a6)} \, \Leftrightarrow \, \mbox{(a1)}.
\]
See \cite[Chapter~3]{LTYen05} for the equivalences between (a5), (a6), and (a1); indeed the equivalence of the former two conditions
is through (a1); see \cite[Proposition~3.3.7]{FacchineiPang03}
regarding the connection between copositivty in (a2) and isolated stationarity which yields the implication
(a2) $\Rightarrow$ (a3).   \hfill $\Box$

Theorems 1 in \cite{Borwein82} offers an extension of the equivalence of the conditions (a6) and (a1) to a non-polyhedral convex feasible
set $P$.  For the implication (a6) $\Rightarrow$ (a1) to be valid in this extended case, the tangent cone ${\cal T}(\bar{x};P)$
in the definition of the critical cone is replaced by the smaller feasible cone of $P$ at $\bar{x}$.

\section{Algebraic Representation of PLQ Functions} \label{sec:representation PLQ}

In \cite{Sun86,Sun92}, Sun explored the structure of convex PLQ functions and obtained a number of fundamental structural results.  Apart from these
early papers, the treatise \cite{RockafellarWets09} has extensive discussion exploring variational properties of PLQ functions; see for instance
Proposition~12.30 and Example~12.31 in the latter reference and also \cite[Theorem~6.1]{BurkeEngle18}.  Our goal in
this section is different: we plan to examine the extension of the max-min representation (\ref{eq:Scholtes PA}) of PA functions to (not necessarily convex) PLQ
functions by starting
with the difference-of-convex representation of piecewise functions in Proposition~\ref{pr:Maher PQ}.  We are also motivated by the elementary representation
of a quadratic function as sums and differences of squared affine functions plus an affine function.
Namely, for a symmetric matrix $Q = P^TDP$ where $D$ is a diagonal matrix
with diagonal entries $\{ d_i \}_{i=1}^n$ and $P$ is an $n \times n$ matrix with rows $\{ P_{i \bullet} \}_{i=1}^n$, then
\begin{equation} \label{eq:quadratic decomposition}
q(x) \, = \, \thalf \, x^TQx + c^Tx + \alpha \, = \, \thalf \, \left[ \, \displaystyle{
\sum_{i \, : \, d_i \, > 0 \, }
} \, d_i \, \left( \, P_{i \bullet} x \, \right)^2 - \displaystyle{
\sum_{i \, : \, d_i \, < 0 \, }
} \, | \, d_i \, | \, \left( \, P_{i \bullet} x \, \right)^2 \, \right] + c^Tx + \alpha
\end{equation}
expresses the quadratic function $q(x)$ as described.  Thus, quadratic functions are composed of two simple classes of convex functions: squares of
linear functions and affine functions, combined together by addition and subtraction.  Using solely squares of affine functions as the ``building functions''
and relying on addition and subtraction only are not enough to yield all PLQ funtions.  The simple squared plus function $t_+^{\, 2}$ suggests that
we need to expand the affine functions to include the most basic PA function, i.e., the plus function; further,
the representation (\ref{eq:Scholtes PA}) suggests that we need to include the pointwise max-min operation.

{\bf The elementary building functions}.  We call the composition of the plus function with an affine function a Plus-Composite-Affine (or in short, PCA)
function; this is a function of the form
$\max(a^Tx + \alpha,0)$ for some vector $a$ and scalar $b$.  One immediate difference between the family of PCA functions and the
family of affine functions is that the latter family is closed under addition and subtraction whereas the former is not.  For our purpose,
we are also interested in the squared PCA functions.  Let ${\cal F}$ consist of two families of functions: squares of affine functions and
squares of PCA functions.  Each member function in ${\cal F}$  is nonnegative, convex, and differentiable.

We begin with a lemma about the distance function to a closed set.  Worthy of note about this lemma is that we employ a polyhedral norm to define
the distance function.  To be specific, we employ the 1-norm:  for a closed set $S \subseteq \mathbb{R}^n$, let
\[
\mbox{dist}_1(x;S) \, \triangleq \, \displaystyle{
\operatornamewithlimits{\mbox{minimum}}_{s \in S}
} \ \| \, s - x \, \|_1.
\]
We should note that a result \cite[Proposition~12.31 part (c)]{RockafellarWets09} related to the one below employs the squared Euclidean-norm distance
function to characterize a polyhedral set.  The lemma characterizes a piecewise polyhedral set in terms of the 1-norm distance function defined above.

\begin{lemma} \label{lm:piecewise polyhedral and PA} \rm
A closed set $S \subseteq \mathbb{R}^n$ is piecewise polyhedral if and only if $\mbox{dist}_1(x;S)$ is a piecewise affine function on $\mathbb{R}^n$.
\end{lemma}
{\bf Proof.}  ``Only if''.  In general, if a closed set $S$ is the union of finitely many closed sets $\{ S^{\, i} \}_{i=1}^I$, then
$\mbox{dist}_1(x,S) = \displaystyle{
\min_{1 \leq i \leq I}
} \, \mbox{dist}_1(x,S^{\, i})$.  Thus the ``only if'' statement follows readily because the 1-norm distance function to a polyhedron is the value
function of a parametric linear program, thus is piecewise affine by well-known linear programming theory.

``If''.  By the max-min representation (\ref{eq:Scholtes PA}) and the nonnegativity of the distance function, it follows that there exist 
affine functions
$\{ f_{ij}(x) \}_{j=1}^{J_i}$ for some positive integer $J_i$ and for all $i = 1, \cdots, I$ for some positive integer $I$ such that
\[
\mbox{dist}_1(x;S) \, = \, \displaystyle{
\max_{1 \leq i \leq I}
} \, \displaystyle{
\min_{1 \leq j \leq J_i}
} \, f_{ij}(x)_+, \epc x \, \in \, \mathbb{R}^n.
\]
Since $S$ is the zero set of the distance function, we deduce that
\[
S \, = \, \displaystyle{
\bigcap_{1 \leq i \leq I}
} \, \left\{ \, x \, \in \, \mathbb{R}^n \, \mid \, \displaystyle{
\min_{1 \leq j \leq J_i}
} \, f_{ij}(x)_+ \, = \, 0 \, \right\} \, = \,  \displaystyle{
\bigcap_{1 \leq i \leq I}
} \, \displaystyle{
\bigcup_{1 \leq j \leq J_i}
} \, \underbrace{\left\{ \, x \, \in \, \mathbb{R}^n \, \mid \, f_{ij}(x) \, \leq \, 0 \, \right\}}_{\mbox{denoted $S^{\, ij}$}},
\]
by the nonnegativity of $f_{ij}(x)_+$.
Since each $S^{\, ij}$ is a halfplane, it follows readily that $S$ is the union of finitely many polyhedra.  \hfill $\Box$

In the following result, we keep the quadratic functions that
define the pieces of the PLQ function in its representation; each such quadratic function has the elementary
decomposition (\ref{eq:quadratic decomposition})
into sums and differences of squared affine functions plus an affine function that can be employed in
(\ref{eq:PLQ representation}) to refine this decomposition.

\begin{proposition} \label{pr:PLQ} \rm
Let $f : \mbox{dom}\, f \to \mathbb{R}$ be a PLQ function on a polyhedral $\mbox{dom}\, f$ that is the union of
finitely many polyhera $\{ P^{\, i} \}_{i=1}^I$; on each such polyhedral piece $P^{\, i}$ is a quadratic function
$q_i$ such that $f(x) = q_i(x)$ for all $x \in P^{\, i}$.  Then there exists finitely many functions $\{ f_{i \wh{k}} \}_{\wh{k}=1}^{\wh{K}_i}$
for $i = 1, \cdots, I$, each given by
\[
f_{i\wh{k}}(x) \, = \, \displaystyle{
\sum_{j=1}^{J_{i \wh{k}}^+}
} \, f_{i\wh{k}j}^+(x) - \displaystyle{
\sum_{j=1}^{J_{i\wh{k}}^-}
} \, f_{i\wh{k}j}^-(x), \epc \wh{k} \, = \, \cdots, \wh{K}_i,
\]
where each $f_{i\wh{k}j}^{\pm} \in {\cal F}$ such that
\begin{equation} \label{eq:PLQ representation}
f(x) \, = \, \displaystyle{
\min_{1 \leq i \leq I}
} \, \left[ \, q_i(x) + \displaystyle{
\max_{1 \leq \wh{k} \leq \wh{K}_i}
} \, f_{i\wh{k}}(x)\, \right] \epc \mbox{for all } x \, \in \, \mbox{dom}\, f,
\end{equation}
and the zero set of the function $\wh{\phi}_i(x) \triangleq \displaystyle{
\max_{1 \leq \wh{k} \leq \wh{K}_i}
} \, f_{i\wh{k}}(x)$ coincides with $P^{\, i}$.
\end{proposition}
{\bf Proof.}  We first remark that the 2-norm in (\ref{eq:Maher PQ cd representation}) can be replaced by the 1-norm; this
replacement results in the following representation of $f(x)$ for all $x \in \mbox{dom}\, f$,
\begin{equation} \label{eq:PQ 1-norm representation}
f(x) \, = \, \displaystyle{
\min_{1 \leq i \leq I}
} \, \left\{ \, q_i(x) + \mbox{dist}_1(x;P^{\, i}) \, \displaystyle{
\max_{1 \leq j \leq I}
} \, \| \nabla q_j(x) - \nabla q_i(x) \|_1 + \displaystyle{
\frac{3 \, L_i}{2}
} \, \left[ \, \mbox{dist}_1(x;P^{\, i}) \, \right]^2 \, \right\} .
\end{equation}
The proof of this identity follows that of (\ref{eq:Maher PQ cd representation}).  In fact, with
\[
\phi_{2i}(x) \, \triangleq \, \mbox{dist}_2(x;P^{\, i}) \, \displaystyle{
\max_{1 \leq j \leq I}
} \, \| \nabla q_j(x) - \nabla q_i(x) \|_2 + \displaystyle{
\frac{3 \, L_i}{2}
} \, \left[ \, \mbox{dist}_2(x;P^{\, i}) \, \right]^2,
\]
the proof of (\ref{eq:Maher PQ cd representation}) hinges on
two things: $\phi_{2i}(x) = 0$ if and only if $x \in P^{\, i}$, and $q_i(x) + \phi_{2i}(x) \geq f(x)$ for all
$x \in \mbox{dom}\, f \setminus P^{\, i}$.
Clearly, these two properties of the functions $\phi_{2i}(x)$ remain valid if we replace them by:
\[
\wh{\phi}_i(x) \, \triangleq \, \mbox{dist}_1(x;P^{\, i}) \, \displaystyle{
\max_{1 \leq j \leq I}
} \, \| \nabla q_j(x) - \nabla q_i(x) \|_1 + \displaystyle{
\frac{3 \, L_i}{2}
} \, \left[ \, \mbox{dist}_1(x;P^{\, i}) \, \right]^2,
\]
because $\| a \|_1 \geq \| a \|_2$ for any vector $a \in \mathbb{R}^n$.  Hence we obtain the 1-norm representation
(\ref{eq:PQ 1-norm representation}) of $f$.  The advantage of the latter representation over the former one is
that we have
\[
\wh{\phi}_i(x) \, = \, \displaystyle{
\max_{1 \leq j \leq I}
} \, \underbrace{\mbox{dist}_1(x;P^{\, i})}_{\mbox{denoted $f_i(x)$}} \, \left[ \,
\underbrace{\| \nabla q_j(x) - \nabla q_i(x) \|_1 + \displaystyle{
\frac{3 \, L_i}{2}
} \, \mbox{dist}_1(x;P^{\, i})}_{\mbox{each denoted $g_{ji}(x)$}} \, \right],
\]
which is the pointwise maximum of finitely many products each of two nonnegative, convex, PA functions.
Next, we examine each such product $f_i(x) g_{ji}(x)$.
By \cite[Theorem~2.49]{RockafellarWets09}, we can write
\[
f_i(x) \, = \, \displaystyle{
\max_{1 \leq k \leq K_i}
} \, \max( \, \ell_{ki}(x), \, 0 \, ) \epc \mbox{and} \epc
g_{ji}(x) \, = \, \displaystyle{
\max_{1 \leq k \leq K_{ji}}
} \, \max( \, \ell_{kji}(x), \, 0 \, ),
\]
where $\ell_{ki}(x)$ and $\ell_{kji}(x)$ are affine functions.  We have
\[ \begin{array}{lll}
f_i(x) \, g_{ji}(x) \, = \,  \displaystyle{
\max_{1 \leq k_i \leq K_i}
} \,  \displaystyle{
\max_{1 \leq k_{ji} \leq K_{ji}}
} \, \left[ \max( \, \ell_{k_ii}(x), \, 0 \, ) \, \max( \, \ell_{k_{ji}ji}(x), \, 0 \, ) \, \right] \\ [0.15in]
\epc = \, \thalf \, \displaystyle{
\max_{1 \leq k_i \leq K_i}
} \,  \displaystyle{
\max_{1 \leq k_{ji} \leq K_{ji}}
} \, \left[ \begin{array}{l}
\left[ \, \max( \, \ell_{k_ii}(x), \, 0 \, ) + \max( \, \ell_{k_{ji}ji}(x), \, 0 \, ) \, \right]^2 \\ [0.1in]
- \, \left[ \, \max( \, \ell_{k_ii}(x), \, 0 \, ) \, \right]^2 - \left[ \, \max( \, \ell_{k_{ji}ji}(x), \, 0 \, ) \, \right]^2
\end{array} \right] \\ [0.3in]
\epc = \, \thalf \, \displaystyle{
\max_{1 \leq k_i \leq K_i}
} \,  \displaystyle{
\max_{1 \leq k_{ji} \leq K_{ji}}
} \, \left[ \begin{array}{l}
\max\left\{ \, \left[ \, \max( \, \ell_{k_ii}(x) + \ell_{k_{ji}ji}(x), \, 0 \, ) \, \right]^2, \ \ell_{k_ii}(x)^2, \ \ell_{k_{ji}ji}(x)^2 \, \right\} \\ [0.2in]
\epc - \, \left[ \, \max( \, \ell_{k_ii}(x), \, 0 \, ) \, \right]^2 - \left[ \, \max( \, \ell_{k_{ji}ji}(x), \, 0 \, ) \, \right]^2
\end{array} \right].
\end{array}
\]
Since for any scalar $t$, we have $t^2 = \max( t,0 )^2 + \max( -t,0 )^2$, we deduce that
$\wh{\phi}_i(x)$ is equal to the pointwise maximum function:
\[ \begin{array}{l}
\displaystyle{
\max_{1 \leq j \leq I}
} \, \displaystyle{
\max_{1 \leq k_i \leq K_i}
} \,  \displaystyle{
\max_{1 \leq k_{ji} \leq K_{ji}}
} \, \left[ \begin{array}{l}
\max\left\{ \, \left[ \, \max( \, \ell_{k_ii}(x) + \ell_{k_{ji}ji}(x), \, 0 \, ) \, \right]^2, \ \ell_{k_ii}(x)^2, \ \ell_{k_{ji}ji}(x)^2 \, \right\} \\ [0.1in]
\epc - \, \left[ \, \max( \, \ell_{k_ii}(x), \, 0 \, ) \, \right]^2 - \left[ \, \max( \, \ell_{k_{ji}ji}(x), \, 0 \, ) \, \right]^2
\end{array} \right] \\ [0.3in]
\epc = \, \displaystyle{
\max_{1 \leq j \leq I}
} \, \displaystyle{
\max_{1 \leq k_i \leq K_i}
} \,  \displaystyle{
\max_{1 \leq k_{ji} \leq K_{ji}}
} \, \max\left\{ \begin{array}{l}
\wh{f}_{ij k_i k_{ji}}(x), \\ [0.1in]
\left[ \, \max( \, -\ell_{k_ii}(x), \, 0 \, ) \, \right]^2 - \left[ \, \max( \, \ell_{k_{ji}ji}(x), \, 0 \, ) \, \right]^2, \\ [0.1in]
\left[ \, \max( \,- \ell_{k_{ji}ji}(x), \, 0 \, ) \, \right]^2 - \left[ \, \max( \, \ell_{k_ii}(x), \, 0 \, ) \, \right]^2
\end{array} \, \right\},
\end{array} \]
where
\[
\wh{f}_{ij k_i k_{ji}}(x) \, \triangleq \,
\left[ \, \max\left( \, \ell_{k_ii}(x) + \ell_{k_{ji}ji}(x), \, 0 \, \right) \, \right]^2 - \left[ \, \max( \, \ell_{k_ii}(x), \, 0 \, ) \, \right]^2
- \left[ \, \max( \, \ell_{k_{ji}ji}(x), \, 0 \, ) \, \right]^2
\]
from which the claimed representation (\ref{eq:PLQ representation}) follows readily.  \hfill $\Box$

\begin{remark} \rm
The above proof provides the following necessary and sufficient representation of a PLQ function.  Namely,
a function $f : \mathbb{R}^n \to \mathbb{R}$ is a PLQ function on a piecewise polyhedral $\mbox{dom}\, f$
if and only if there exist a family of
quadratic functions $\{ q_i \}_{i=1}^I$ and two families of piecewise affine functions $\{ \wh{h}_i(x) \}_{i=1}^I$
and $\{ \, \wt{h}_i(x) \, \}_{i=1}^I$ such that $\mbox{dom}\, f \, \subseteq \,
\displaystyle{
\bigcup_{i=1}^I
} \, \{ x \in \mathbb{R}^n \, \mid \, f(x) = q_i(x) \}$;

\[
f(x) \, = \, \displaystyle{
\min_{1 \leq i \leq I}
} \, \left[ \, q_i(x) + \wh{h}_i(x) \, \wt{h}_i(x) \, \right], \epc \forall \, x \, \in \, \mbox{dom}\, f
\]
and $\{ x \in \mbox{dom}\, f \, \mid \, f(x) \, = \, q_i(x) \} = \{ x \in \mbox{dom}\, f \, \mid \, \wh{h}_i(x) \, = \,0 \}$
for each $i = 1, \cdots, I$. \hfill $\Box$
\end{remark}

\section{Second-Order Properties of  Piecewise Quadratic Functions} \label{sec:PQ functions}

In this session, we discuss the second-order directional properties of PC$^{\, 2}$ functions.  The results herein are not
surprising and yet seemingly new.

\begin{proposition} \label{pr:PC are dd2} \rm
Let $f$ be a PC$^{\rm 2}$ function on an open set $\Omega \subseteq \mathbb{R}^n$.  Then $f$ is twice directionally differentiable on $\Omega$.
Moreover, for every pair $(x,d) \in \Omega \times \mathbb{R}^n$, $f^{(2)}(x;d)$ is equal to $d^{\, T} \nabla ^2f_i(x)d$ for
any $i \in {\cal A}^{\, \prime}(x;d)$.
\end{proposition}
{\bf Proof.}  The proof follows the line of proof of Lemma~4.6.1 in \cite{FacchineiPang03} cited above.
As in this lemma, it suffices to show that $d^{\, T} \nabla^2 f_i(x)d = d^{\, T} \nabla^2 f_j(x)d$ for any two
indices $i$ and $j$ in ${\cal A}^{\, \prime}(x;d)$.  Assume the contrary.  Let $\bar{i}$ and $\bar{j}$ be two indices
in  ${\cal A}^{\, \prime}(x;d)$ such that $d^{\, T} \nabla^2 f_{\bar{i}}(x)d \neq d^{\, T} \nabla^2 f_{\bar{j}}(x)d$.  Since
$f_{\bar{i}}$ and $f_{\bar{j}}$ are C$^{\, 2}$ functions,
$f_{\bar{i}}(x) = f_{\bar{j}}(x)$, and
$\nabla f_{\bar{i}}(x)^Td = \nabla f_{\bar{j}}(x)^Td$, it follows that a scalar $\varepsilon_{\bar{i} \bar{j}} > 0$
exists such that  $f_{\bar{i}}(x + \tau d) \neq f_{\bar{j}}(x + \tau d)$ for all $\tau \in ( \, 0, \varepsilon_{\bar{i} \bar{j}} \, ]$.
At this point, the same proof of Lemma~4.6.1 in \cite{FacchineiPang03} can be applied to derive a contradiction; in essence,
this argument relies solely on the compactness of the line segment $\left[ x, x + \varepsilon \, d \, \right]$, where
$\varepsilon > 0$ is a suitable scalar derived from the $\varepsilon_{\bar{i} \bar{j}}$, appropriately reduced if necessary to
ensure that ${\cal A}(x + \tau d) \subseteq {\cal A}(x)$ for all $\tau \in [ 0, \varepsilon ]$.  We omit the details.
\hfill $\Box$

\begin{remark} \label{re:dd and semi} \rm
Although the deficiencies of the second directional derivative
$f^{(2)}(x;d)$ have been very well noted in \cite[Section~13.B]{RockafellarWets09},  Proposition~\ref{pr:PC are dd2} suggests
that twice directional differentiability is a weaker requirement than twice
semidifferentiability in that the derivative $f^{(2)}(x;d)$ may exist while the second-order limit (\ref{eq:2-order semidiff limit})
does not.    As asserted by Proposition~\ref{pr:PC are dd2}, a PC$^2$ function is always twice directionally differentiable;
but it may not be twice semidifferentiable.  One counterexample is given in \cite[Example 13.10]{RockafellarWets09},
where $f(x) = \max(|x+a|^2,1)$ with $|a| = 1$ is a univariate PQ function that fails to be twice semidifferentiable at $x = 0$.
An example at the end of this section further illustrates the difference between these two second-order differentiability
concepts.   \hfill $\Box$
\end{remark}

It is interesting to compare Propositions~\ref{pr:PC are dd2} with \ref{prop:convex plq}.
In the latter proposition (for PLQ functions), we obtained the second directional derivative $f^{(2)}(\bar{x};v)$ for all $v \in {\cal T}(\bar{x};P^{\, i})$, whereas
in the former proposition (for PC$^{\, 2}$ functions), it is not difficult to see that $i \in {\cal A}^{\, \prime}(x;d)$ if and only if
$d \in {\cal R}(x;P^{\, i})$, which is the so-called ``radial cone'' of the (not necessarily polyhedral) piece $P^{\, i}$ that is a subset of
the tangent cone ${\cal T}(x;P^{\, i})$.  Thus Proposition~\ref{pr:PC are dd2} gives the second directional derivative $f^{(2)}(x;d)$
for all $d \in {\cal R}(x;P^{\, i})$.  The two cones ${\cal R}(x;P^{\, i})$ and ${\cal T}(x;P^{\, i})$ coincide when $P^{\, i}$ is polyhedral.
If $P^{\, i}$ is convex for $i \in {\cal A}(x)$, then $P^{\, i} \subseteq x + {\cal R}(x;P^{\, i})$.

The next proposition generalizes the result of Proposition \ref{prop:convex plq} on the local exactness of the quadratic expansion of a
PQ function restricted to directions in the radial cones at a point.

\begin{proposition} \label{pr:PLQ} \rm  Let $f$ be a PQ function on a domain ${\cal D} \subseteq \mathbb{R}^n$.
Then, for every $\bar{x} \in {\cal D}$ and every piece $P^{\, i}$ of $f$ containing $\bar{x}$,
it holds that for all $x \in \bar{x} + {\cal R}(\bar{x};P^{\, i})$,
\begin{align}\label{eq:PLQequality}
 f(x) \, = \, f(\bar{x}) + f^{\, \prime}(\bar{x}; x - \bar{x}) + \thalf \, f^{(2)}(\bar{x};x - \bar{x}),
\end{align}
Thus $f^{(2)}(\bar{x};\bullet)$ is continuous when restricted to the cone ${\cal R}(\bar{x};P^{\, i})$.
\end{proposition}
{\bf Proof.}
If $x \in \bar{x} + {\cal R}(\bar{x};P^{\, i})$, then
$i \in {\cal A}^{\, \prime}(\bar x; x-\bar x)$.  This implies that $f^{\, \prime}(\bar x; x-\bar x) = \nabla f_i(\bar x)^T(\bar{x}; x - \bar{x})$
Lemma~4.6.1 in \cite{FacchineiPang03} and that $f^{(2)}(x;x-\bar x)=(x - \bar{x})^T \nabla^2f_i(\bar x)(\bar{x};x - \bar{x})$
by Proposition \ref{pr:PC are dd2}.  Since for the quadratic function $f_i$, we have
\begin{align}
f_i(x) \, = \, f_i(\bar{x}) + \nabla f_i(\bar x)(\bar{x}; x - \bar{x}) + \thalf \,(x - \bar{x})^T \nabla^2f_i(\bar x)(\bar{x};x - \bar{x}),
\end{align}
(\ref{eq:PLQequality}) follows readily.  The last statement of the proposition is obvious.  \hfill $\Box$

\begin{remark} \rm
If $P^{\, i}$ is convex, then (\ref{eq:PLQequality}) holds for all $x \in P^{\, i}$.  Thus, if $f$ is PQ with convex pieces, then (\ref{eq:PLQequality})
holds for all $x$ near $\bar{x}$.  This raises a question that we will formally pose in the next subsection and for which we do not have an answer presently.
Nevertheless, the next proposition gives a partial answer to this question. \hfill $\Box$
\end{remark}

If a PQ function is continuously differentiable, then it is also a PLQ function.  This seems to be a new result in the literature of PQ functions.

\begin{proposition} \label{pr:PQ2} \rm
Let $f :\Omega \to \mathbb{R}$ be a C$^{\, 1}$ function defined on the open set $\Omega$ containing $\bar{x}$.  The following
three statements are equivalent.\\[0.05in]
(a) $f$ is piecewise quadratic near $\bar{x}$;\\[0.05in]
(b) $\nabla f$ is a piecewise affine near $\bar{x}$;\\[0.05in]
(c) $f$ is piecewise linear-quadratic near $\bar{x}$.

\end{proposition}
{\bf Proof.}  (a) $\Rightarrow $ (b).  This follows from \cite[Lemma~2]{Rockafellar03}.

(b) $\Rightarrow$ (c).  Write $F(x) \triangleq \nabla f(x)$.  Let $\{ A^i x + b^i \}_{i=1}^I$
and $\{ P^{\, i} \}_{i=1}^I$ be the affine pieces of $F$ in a neighborhood ${\cal N}$ of $\bar{x}$ that we may assume to be polyhedral
such that $F(x) = A^ix+ b^i$ for all $x \in {\cal N} \cap P^{\, i}$, where each $A^i \in \mathbb{R}^{n \times n}$, $b^i \in \mathbb{R}^n$
and $P^{\, i}$ is a polyhedral set.  We may assume without loss of generality that this neighborhood ${\cal N}$ is such that
\[
F(x) \, = \, F(\bar{x}) + F^{\, \prime}(\bar{x};x - \bar{x}), \epc \forall \, x \, \in \, {\cal N}.
\]
Since PA functions are semismooth \cite[Definition~7.4.2]{FacchineiPang03}, it follows that $f$ is SC$^{\, 1}$
at $\bar{x}$ \cite[Section~7.4.1]{FacchineiPang03}.
From expression~(7.4.14) in \cite{FacchineiPang03} for a SC$^{\, 1}$ function, we deduce
\[
\displaystyle{
\lim_{\tau \downarrow 0}
} \, \displaystyle{
\frac{f(\bar{x} + \tau v) - f(\bar{x}) - \tau \, \nabla f(\bar{x})^Tv - \displaystyle{
\frac{\tau^2}{2}
} \, v^TF^{\, \prime}(\bar{x};v)}{\tau^2}
} \, = \, 0,
\]
which readily yields that $f^{(2)}(\bar{x};v) = v^T F^{\, \prime}(\bar{x};v)$
for all $v \in \mathbb{R}^n$.   Since $F^{\, \prime}(\bar{x};\bullet)$ is a PA function on $\mathbb{R}^n$, it
follows from \cite[Proposition~4.2.1]{FacchineiPang03} that there exists a ``polyhedral subdivision'' $\Xi$ of
$\mathbb{R}^n$ such that $F^{\, \prime}(\bar{x};\bullet)$ coincides with one of the linear function $\{ A^iv \}_{i=1}^I$
on each polyhedron in $\Xi$.  
Letting $\{ \wh{P}^{\, j} \}_{j=1}^J$ be the polyhedra in the subdivision $\Xi$, we deduce that
$f^{(2)}(\bar{x};\bullet)$ is a (homogenous) quadratic function on each $\wh{P}^{\, j}$.  More precisely, for each $j = 1, \cdots, J$,
there exists $i_j \in \{ 1, \cdots, I \}$ such that $f^{(2)}(\bar{x};v) = v^TA^{i_j}v$ for all $v \in \wh{P}^{\, j}$.
By showing that (\ref{eq:PLQequality}) holds for all $x$ in ${\cal N}$, it will imply
that $f$ is piecewise linear-quadratic near $\bar{x}$.
For a fixed but arbitrary $x \in {\cal N}$, define the univariate function $\psi(t) \triangleq f(\bar{x} + t (x - \bar{x})) - f(\bar{x}) - t \nabla f(\bar{x})^T( x - \bar{x})$
for $t \in [ 0,1 ]$.  This function is differentiable with derivative
\[
\psi^{\, \prime}(t) \, = \, \left[ \, F(\bar{x} + t(x - \bar{x})) - F(\bar{x}) \, \right]^T( x - \bar{x}) \, = \, t \, F^{\, \prime}(\bar{x};x - \bar{x})^T( \, x - \bar{x} \, ).
\]
Hence,
\[ \begin{array}{lll}
f(x) - f(\bar{x}) -  \nabla f(\bar{x})^T( x - \bar{x}) & = & \psi(1) - \psi(0) \\ [5pt]
& = & \displaystyle{
\int_0^1
} \, \psi^{\, \prime}(t) \, dt \, = \, \thalf \, F^{\, \prime}(\bar{x};x - \bar{x})^T ( \, x - \bar{x} \, ),
\end{array} \]
which is the desired equality (\ref{eq:PLQequality}).

(c) $\Rightarrow$ (a).  This is obvious.  \hfill $\Box$

The example below illustrates many of the results establish above.

\begin{example} \label{ex:PQ example} \rm
Consider the following piecewise quadratic function:
\begin{equation} \label{eq:PQ example}
f(x) \, = \, \thalf \, \left[ \, \max\left( \, \| \, x \, \|_2^2, \, 1 \, \right) - x^TQx \, \right], \epc x \, \in \, \mathbb{R}^n,
\end{equation}
where $Q$ is a symmetric matrix, which is not necessarily positive semidefinite.  One piece of this function is the exterior of the
unit ball, thus not convex.  It is not difficult to verify
the following directional derivatives of the first and second order: for every pair $(x,d) \in \mathbb{R}^n \times \mathbb{R}^n$,
\begin{equation} \label{eq:derivatives in example}
\begin{array}{l}
f^{\, \prime}(x;d) \, = \, \left\{ \begin{array}{ll}
x^Td - x^TQd & \; \mbox{if $\| \, x \, \|_2 > 1$} \\ [5pt]
-x^TQd & \; \mbox{if $\| \, x \, \|_2 < 1$} \\ [5pt]
\max\left( \, x^Td, 0 \, \right) - x^TQd & \; \mbox{if $\| \, x \, \|_2 = 1$};
\end{array} \right. \\ [0.4in]
{\rm d}^2 f(x)(d) \, = \, f^{\rm (2)}(x;d) \, = \, \left\{ \begin{array}{ll}
\| \, d \, \|_2^2 - d^{\, T}Qd & \; \mbox{if $\| \, x \, \|_2 > 1$ or [ $\| \, x \, \|_2 = 1$ and  $x^Td \, \geq \, 0$ ]} \\ [5pt]
-d^{\, T}Qd & \; \mbox{if $\| \, x \, \|_2 < 1$ or [ $\| \, x \, \|_2 = 1$ and  $x^Td \, < \, 0$ ]}.
\end{array} \right.
\end{array} \end{equation}
Both second-order directional derivatives $f^{\rm (2)}(x;d) = {\rm d}^2 f(x)(d)$ exist for all $(x,d)$
and yet are discontinuous in neither variable while the other is fixed.  Thus this PQ function $f$
is  not twice semidifferentiable.  \hfill $\Box$
\end{example}

\subsection{Some open questions} \label{subsec:open questions}

The results in this section and Section~\ref{sec:representation PLQ} have added to the understanding of PLQ and PQ functions.
Yet, there remain several questions whose answers we do not know at this time and which seem worthwhile to ask for future
research.  The main question is whether we can characterize a PQ function to be PLQ in terms of several properties of the latter.
The following are some specific questions:\\[0.05in]
$\bullet $ If the domain of a PQ function is the union of finitely many closed {\sl convex} sets on each of which the function is
quadratic, does it follow that the PQ function is PLQ?\\[0.05in]
$\bullet $ If a PQ functions is twice semidifferentiable, is it necessarily a PLQ function?\\[0.05in]
$\bullet $ Is there a ``simpler'' representation of a PLQ function in terms of the family of functions in ${\cal F}$ introduced prior to Proposition~\ref{pr:PLQ}
than the one (\ref{eq:PLQ representation}) in this proposition?\\[0.05in]
$\bullet $ Is the class of functions with the representation (\ref{eq:PLQ representation}) equal to the class of PQ functions?

\section{Second-Order Optimality Conditions}

Our goal in this section is to extend the optimality results in Subsection~\ref{subsec:QP} to a linearly constrained piecewise linear-quadratic program.
For simplicity, in both Theorems~\ref{th:qp SONC} and \ref{th:qp SOSC}, we take the objective $f$ to be a PLQ function on the entire $\mathbb{R}^n$.  As such,
$f$ is twice semidifferentiable on $\mathbb{R}^n$.
The first result concerns a local minimizer that extends Proposition~\ref{pr:QP local}.

\begin{theorem} \label{th:qp SONC} \rm
Let $f : \mathbb{R}^n \to \mathbb{R}$ be a PLQ function with polyhedral pieces $\{ P^{\, i} \}_{i=1}^I$
and associated quadratic functions $\{ \, q_i \, \}_{i=1}^I$.
Let $X$ be a polyehedral set in $\mathbb{R}^n$.  Let $\wt{P}^{\, i} \, \triangleq X \, \cap \, P_i$.
The following four statements are all equivalent at a given vector $\bar{x} \in X$:
\\[0.05in]
(a1) $\bar{x}$ is a local minimizer of $f$ on $X$;
\\[0.05in]
(a2) for every $i \in {\cal A}(\bar{x})$, $\bar{x}$ is a local minimizer of $q_i$ on $\wt{P}^{\, i}$;
\\[0.05in]
(b1) $\bar{x}$ is a d-stationary point of \eqref{eq:basic opt problem} and satisfies the second-order necessary condition; 
\\[0.05in]
(b2) for every $i \in {\cal A}(\bar{x})$, $\bar{x}$ is a stationary point of $f$ (or equivalently, $q_i$) on $\wt{P}^{\, i}$
and $\nabla^2 q_i(\bar{x})$
is copositive on ${\cal C}(\bar{x};q_i;\wt{P}^{\, i})$.
%
\end{theorem}
{\bf Proof.} 
%
(a1) $\Rightarrow$ (a2): Let $\mathcal{N}$ be a neighborhood of $\bar{x}$ such that ${\cal A}(x) \subseteq {\cal A}(\bar{x})$
for all $x \in {\cal N}$.
We claim that for any $i \in {\cal A}(\bar{x})$, $\bar{x}$ is a minimizer of $f$ on ${\cal N} \cap \wt{P}^{\, i}$.
Indeed, for any such $i$, we have $q_i(x) = f(x) \geq f(\bar{x}) = q_i(\bar{x})$ for any
$x\in \mathcal{N}\cap \widetilde{P}^{\,i}$.

(a2) $\Rightarrow$ (a1):  Choose a neighborhood ${\cal N}$ of $\bar{x}$ satisfying two conditions: (i) $\bar{x}$ is a
minimizer of $q_i$ on ${\cal N} \cap \widetilde{P}^{\,i}$ for every $i \in {\cal A}(\bar{x})$, and (ii)
${\cal A}(x) \subseteq {\cal A}(\bar{x})$
for every $x \in {\cal N}$.  Let $x \in X \cap {\cal N}$ be arbitrary.  For every $i \in {\cal A}(x)$, we have
\[
f(x) \, = \, q_i(x) \, \geq \, q_i(\bar{x}) \, = \, f(\bar{x}),
\]
where the equalities hold by the choice of $i$ and the local minimizing property of $\bar{x}$ for $q_i$ on each piece $\widetilde{P}^{\,i}$.

(b1) $\Rightarrow$ (b2):  This holds because $v\in \mathcal{T}(\bar{x}; \widetilde{P}^{\,i})$ for some
$i\in \mathcal{A}(\bar{x})$ implies $v\in \mathcal{T}(\bar{x}; X)\cap \mathcal{T}(\bar{x}; {P}^{\,i})$, which, by Proposition \ref{prop:convex plq}, further yield
\begin{equation}\label{proof:eq1}
f^{\, \prime}(\bar{x};v) \, = \, \nabla q_i(\bar{x})^Tv\epc \mbox{and}\epc
f^{(2)}(x;v) \, = \, v^T \nabla^2 q_i(\bar x)v.
\end{equation}
(b2) $\Rightarrow$ (b1): This holds because
for any $v\in \mathcal{T}(\bar{x};X)$, if $v\in \mathcal{T}(\bar{x};\widetilde{P}^{\,i})$ for some $i\in \mathcal{A}(\bar{x})$,
then $v\in  \mathcal{T}(\bar{x};P^{\,i})$ and thus \eqref{proof:eq1} holds.

%

(a2) $\Leftrightarrow$ (b2): by Proposition~\ref{pr:QP local}.
\hfill $\Box$

\begin{remark} \rm
While the proof is not difficult, the implication (b1) $\Rightarrow$ (a1) is missing in the literature till now.  Thus
Theorem~\ref{th:qp SONC} gives a complete set of necessary and sufficient conditions for the local optimality of PLQ programs
in terms of the second-order necessary conditions and the copositivity condition (b2).  \hfill $\Box$
\end{remark}

Employing \cite[Theorem~3]{Borwein82}, we can deduce that Theorem~\ref{th:qp SONC} remains valid for a non-polyhedral
constraint set $X$ provided that the tangent cone ${\cal T}(\bar{x};X)$ is polyhedral.  We omit the details.
We next extend Proposition~\ref{pr:QP strong local} to a PLQ program.  The extension relies on the equivalence of the piecewise program
locally to the pieces that contain the point $\bar{x}$ in question, similar to the equivalence of (a1) to (a2) in the above
Proposition~\ref{th:qp SONC}.  Once such a local equivalence is establish, all the other equivalent conditions follow readily
from the previous results for a QP.

%
%

\begin{theorem} \label{th:qp SOSC} \rm
Let $f : \mathbb{R}^n \to \mathbb{R}$ be a PLQ function with polyhedral pieces $\{ P^{\, i} \}_{i=1}^I$
and associated quadratic functions $\{ \, q_i \, \}_{i=1}^I$.
Let $X$ be a polyhedral set in $\mathbb{R}^n$.  Let $\wt{P}^{\, i} \, \triangleq X \, \cap \, P_i$.
The following statements are all equivalent at a given vector $\bar{x} \in X$:\\[0.05in]
(a1) $\bar{x}$ is a strong local minimizer of $f$ on $X$;
\\[0.05in]
(a2) $\bar{x}$ is a strict local minimizer of $f$ on $X$;
\\[0.05in]
(a3) $\bar{x}$ is an isolated local minimizer of $f$ on $X$;
\\[0.05in]
(a4) $\bar{x}$ is an isolated stationary point and a local minimizer of $f$ on $X$;
\\[0.05in]
(b1) for every $i \in {\cal A}(\bar{x})$, $\bar{x}$ is a strong local minimizer of $f$ (or equivalently $q_i$) on $\wt{P}^{\, i}$;
\\[0.05in]
(b2) for every $i \in {\cal A}(\bar{x})$, $\bar{x}$ is a strict local minimizer of $f$ (or equivalently $q_i$) on $\wt{P}^{\, i}$;
\\[0.05in]
(b3) for every $i \in {\cal A}(\bar{x})$, $\bar{x}$ is an isolated local minimizer of $f$ (or equivalently $q_i$) on $\wt{P}^{\, i}$;
\\[0.05in]
(b4) for every $i \in {\cal A}(\bar{x})$, $\bar{x}$ is an isolated stationary point and a local minimizer of $f$ (or equivalently $q_i$) on $\wt{P}^{\, i}$;
\\[0.05in]
(c) ${\rm d}f(\bar{x})(x-\bar{x}) \geq 0$ for all $x\in X$ and ${\rm d}^2 f(\bar{x})(x-\bar{x}) > 0$  for all $x \in X \setminus \{ \bar{x} \}$
with ${\rm d}f(\bar{x})(x-\bar{x}) = 0$;
\\[0.05in]
(d1) $\bar{x}$ is a d-stationary point of \eqref{eq:basic opt problem} and satisfies the second-order sufficient condition; 
\\[0.05in]
(d2) for every $i \in {\cal A}(\bar{x})$, $\bar{x}$ is a stationary point of $q_i$ on $\wt{P}^{\, i}$ and $\nabla^2 q_i(\bar{x})$
is strictly copositive on ${\cal C}(\bar{x};q_i;\wt{P}^{\, i})$.
\end{theorem}
{\bf Proof.}  We may proceed as in the proof of Theorem~\ref{th:qp SONC} to show the equivalence of the individual statements (a1) through (a4) for the
problem (\ref{eq:basic opt problem}) with the corresponding statements (b1) through (b4) for the piecewise programs.  The inter-equivalences among
the statements (b1) through (b4) and their equivalences with (d1) and (d2) are through Proposition~\ref{pr:QP strong local} for a standard QP.  Finally,
the equivalence with (c) is by Proposition~\ref{pr:second-order optimality}. \hfill $\Box$

\begin{remark} \rm
Similar to the previous Theorem~\ref{th:qp SONC}, Theorem~\ref{th:qp SOSC} gives a complete set of necessary and sufficient
conditions for the (strong, strict, isolated) local optimality in a PLQ program in terms of the second-order sufficient conditions and the
strict copositivity condition (d2) on the pieces.
Many implications in Theorem~\ref{th:qp SOSC} remain valid for a PC$^{\, 2}$ function with convex pieces.
Without the PLQ property, however, it is not possible to apply Proposition~\ref{pr:QP strong local} to establish the complete
equivalences; in particular, to show the necessity condition (d2) under either (a1) or (a2).  \hfill $\Box$
\end{remark}


%

\begin{example} \label{eq:counter-example} \rm
We use the function in Example~\ref{ex:PQ example} to illustrates two important points.
\\[0.05in]
$\bullet $ For a piecewise quadratic (as opposed to piecewise linear-quadratic) program, a stationary point
satisfying the $f^{(2)}(x;\bullet)$ (or even ${\rm d}^2f(x)(\bullet)$) based second-order necessary
condition is not necessarily a local minimizer; in other words, for a PQ program, such a
second-order necessary condition is not in general sufficient for local optimality.  Hence the linear-quadratic property of the objective
function is essential for such sufficiency to hold as established in Theorem~\ref{th:qp SONC}.
\\[0.05in]
$\bullet $ The second-order sufficient condition in terms of the second directional derivative $f^{(2)}(x;\bullet)$ or
the second semiderivative $\mbox{d}^2f(x)(\bullet)$ (which are equal for this example) is not sufficient for a local minimizer
when the domain of some piece is not convex.
This confirms that the second-order sufficient condition based on either one of these second derivatives is weaker than that based on
${\rm d}^2 (f + \delta_X)(x|0)(v)$ as established in Theorem~\ref{thm: RW: opt}, the latter offers an elegant yet abstract necessary and sufficient
condition for strong local optimality of a general nonsmooth, nonconvex program without exposing the set $X$.

We first characterize the second-order stationarity conditions based on $f^{(2)}(x;\,\bullet)$.
Let $\bar{x} \in \mathbb{R}^n$ with $\| \bar{x} \|_2 = 1$ be arbitrary.  The following two statements hold for the function
\\[0.05in]
$\bullet $ $\bar{x}$ is an unconstrained (directional) stationary point of $f$ if and only if $\bar{x}$ is a normalized eigenvector
of the matrix $Q$ corresponding to an eigenvalue $\beta \in [ 0,1 ]$;
\\[0.05in]
$\bullet $ if 0 and 1 are not eigenvalues of $Q$, then $\bar{x}$ satisfies the second-order necessary condition
of $f$ if and only if it is stationary and
\begin{equation} \label{eq:2nd order PQ example}
\bar{x}^T d \, = \, 0 \ \Rightarrow \ d^{\, T}Qd \, \leq \, \| \, d \, \|_2^2.
\end{equation}
{\bf Proof.}  By the expression (\ref{eq:derivatives in example}) of $f^{\, \prime}(\bar{x};d)$,
we deduce that $\bar{x}$ is an unconstrained stationary point of $f$ if and only if
\[
\max\left( \, \bar{x}^Td, 0 \, \right) - \bar{x}^TQd \, \geq \, 0, \epc \forall \, d \, \in \, \mathbb{R}^n.
\]
In turn, this is equivalent to two implications:
\[ \begin{array}{lll}
\bar{x}^Td \, \geq \, 0 & \Rightarrow & ( \, \bar{x} - Q \bar{x} \, )^Td \, \geq \, 0 \\ [5pt]
\bar{x}^Td \, \leq \, 0 & \Rightarrow & \bar{x}^TQd \, \leq \, 0.
\end{array} \]
It is not difficult to show that these inequalities are equivalent to the existence of a scalar $\beta \in [ 0,1 ]$ such that $Q\bar{x} = \beta \bar{x}$,
which is equivalent to the claimed eigenvalue characterization of $\bar{x}$.
Further, if $\bar{x}$ is an unconstrained stationary point of $f$
and $d$ is such that $f^{\, \prime}(\bar{x};d) = 0$, then we must have $\bar{x}^Td = 0$.  Hence
if 0 and 1 are not eigenvalues of $Q$, then by the expression of $f^{\rm (2)}(\bar{x};d)$,
it follows that $\bar{x}$ satisfies the second-order necessary condition of $f$ if and only if
$\bar{x}$ is a normalized eigenvector of the matrix $Q$ corresponding to
an eigenvalue $\beta \in ( 0,1 )$ and the implication (\ref{eq:2nd order PQ example}) holds.
\hfill $\Box$

In the rest of the discussion of the example, we let $n = 2$ and $Q$ be a $2 \times 2$ positive diagonal matrix with
diagonal elements $Q_{11}$ and $Q_{22}$ satisfying:
$0 < Q_{22} < Q_{11} < 1$.   We also fix $\bar{x} = ( 0,-1 )$.  Then $\bar{x}$ is a normalized eigenvalue of $Q$ corresponding to $Q_{22}$.
Hence $\bar{x}$ is a directional stationary point of the function $f$ given by (\ref{eq:PQ example}).  Moreover, since the eigenvalues
of $Q$ are both less than unity, it follows that $\bar{x}$ satisfies the second-order necessary condition.
We show however that $\bar{x}$ is not an unconstrained local minimizer of $f$ by considering the points
\[
x(\varepsilon) \, \triangleq \, \left( \, \sqrt{\displaystyle{
\frac{2 \, \varepsilon}{1 + \displaystyle{
\frac{Q_{11}}{Q_{22}}
}}}}, \, -\sqrt{1 - \varepsilon} \, \right), \epc \mbox{for all $\varepsilon > 0$ sufficiently small}.
\]
We have
\[ \begin{array}{lll}
f(x(\varepsilon)) & = & \thalf \, \left[ \, \max\left( \, \displaystyle{
\frac{2 \, \varepsilon}{1 + \displaystyle{
\frac{Q_{11}}{Q_{22}}
}}} + 1 - \varepsilon, \, 1 \, \right) - \left( \,  \displaystyle{
\frac{2 \, \varepsilon}{1 + \displaystyle{
\frac{Q_{11}}{Q_{22}}
}}} \, \right) Q_{11} - ( \, 1 - \varepsilon \, ) \, Q_{22} \, \right] \\ [0.4in]
& = & \thalf \, \left[ \, 1 - Q_{22} - \varepsilon \, \left\{ \, \left( \,  \displaystyle{
\frac{2}{1 + \displaystyle{
\frac{Q_{11}}{Q_{22}}
}}} \, \right) \, Q_{11} - Q_{22} \, \right\} \, \right] \\ [0.4in]
& = &
\thalf \, \left[ \, 1 - Q_{22} - \displaystyle{
\frac{\varepsilon \, Q_{22}}{1 + \displaystyle{
\frac{Q_{11}}{Q_{22}}
}}} \, \left( \, \displaystyle{
\frac{Q_{11}}{Q_{22}}
} - 1 \, \right) \, \right] \, < \, \thalf ( 1 - Q_{22} ) \, = \, f(\bar{x}).
\end{array} \]
Thus, $\bar{x}$ satisfies the second-order necessary condition but is not an unconstrained local minimizer of
the bivariate function $f(x_1,x_2) = \thalf \left[ \max( x_1^2 + x_2^2, \, 1 ) - Q_{11} x_1^2 - Q_{22} x_2^2 \right]$.

Notice that $f^{\, \prime}(\bar{x};d) = 0$ if and only if $d_2 = 0$, then the second order sufficient condition
\[
f^{\rm (2)}(\bar{x};d) = (1-Q_{11})\,d_1^2>0\epc\forall d\neq 0 \; \mbox{such that $d_2=0$},
\]
actually holds at $\bar{x}$. This indicates that the second-order conditions
defined by $f^{(2)}(x;\,\bullet\,)$ may be unfavourable for general PQ programs.

As a comparison, one can derive from the formula of \cite[Example 13.16]{RockafellarWets09} that
\[
{\rm d}^{\,2}f(\bar{x}\,|\,0)(d) = \max_{\lambda} \left\{\,\lambda \|d\|^2 - d^TQ d \,\mid\, \lambda\bar{x} = Q\bar{x}, \, \lambda\in [0,1]\,\right\}.
\]
Then for any $0\neq d\in \mathbb{R}^2$ with $d_2 = 0$, one has
\[
{\rm d}^{\,2}f(\bar{x}\,|\,0)(d)\, =\, \max_{\lambda}\left\{\,\lambda\, d_{\,1}^{\,2} - Q_{11}d_{\,1}^{\,2} \, \mid \,
\lambda = Q_{22} \in (0,1)\,\right\} = (Q_{22} - Q_{11})\,d_1^2.
\]
Since $Q_{22} < Q_{11}$, the second-order necessary condition defined by ${\rm d}^{\,2}f(\bar{x}\,|\,0)(d) \geq 0$ for all $d$ satisfying
$f^{\, \prime}(\bar{x};d) = 0$ fails at $\bar{x}$.  \hfill $\Box$
\end{example}

We give below another easy result that is seemingly new too.  A realization of this result is given by the problem (\ref{eq:log-likelihood})
arising from a log-likelihood piecewise affine estimation problem.

\begin{proposition} \label{pr:cvx with PA} \rm
Let $f = \phi \circ \psi$ be the composite of a convex function $\phi$ and a PA function $\psi$.  With $X$ being a closed convex set, any
(directional) stationary solution of (\ref{eq:basic opt problem}) is a local minimizer.
\end{proposition}
{\bf Proof.}  Let $\bar{x}$ be a (directional) stationary solution of (\ref{eq:basic opt problem}) and $x \in X$ be arbitrary.
We have
\[ \begin{array}{llll}
f(x) & = & \phi(\psi(x)) \\ [5pt]
& \geq & \phi(\psi(\bar{x})) + \phi^{\, \prime}(\psi(\bar{x});\psi(x) - \psi(\bar{x})),
& \mbox{by convexity of $\phi$} \\ [0.1in]
& = &  \phi(\psi(\bar{x})) + \phi^{\, \prime}(\psi(\bar{x});\psi^{\, \prime}(\bar{x};x - \bar{x})), &
\mbox{for all $x$ near $\bar{x}$, by the PA property of $\psi$} \\ [0.1in]
& = & f(\bar{x}) + f^{\, \prime}(\bar{x};x- \bar{x}),
\end{array} \]
where the last equality is by the directional derivative formula of composite functions.  \hfill $\Box$

\subsection{Finite number of strong local minima} \label{subsec:finite}

In this subsection, we establish the interesting result that the number of strong local minima of a quadratic program is finite, from which
the same conclusion holds for a PLQ in view of the equivalence between (a1) and (b1) in Theorem ~\ref{th:qp SOSC} and the fact that
there are only finitely many QP pieces of a PLQ program.  We will subsequently connect the result with an advanced theory of subanalytic functions.

\begin{proposition} \label{pr:finite strong} \rm
For the quadratic program
\begin{equation}\label{opt:qp}
\begin{array}{ll}
\displaystyle\operatornamewithlimits{minimize}_{x\in \mathbb{R}^n} & \thalf \, x^TQx + c^Tx\\[0.05in]
\mbox{subject to}  & Ax \leq b,
\end{array}
\end{equation}
the set of its isolated (equivalently, strict, strong) local minima is finite.
\end{proposition}
{\bf Proof.}
Denote
\[
\mathcal{F} \, \triangleq \, \left\{ \,
\beta \, \subseteq \{1, \ldots, m \} \, \bigg| \,
\begin{array}{l}
\mbox{there exists an isolated local minimizer} \\
\mbox{with a multiplier $\lambda$ such that $ \mbox{supp}(\lambda) = \beta$}
\end{array}
\right\}.
\]
It suffices to show that for any $\bar{\beta}\in \mathcal{F}$, the corresponding isolated local minimizer $\bar{x}$
with a multiplier $\bar{\lambda}$ satisfying $ \mbox{supp}(\bar{\lambda}) = \bar{\beta}$ is unique.
Based on the KKT optimality condition of \eqref{opt:qp} at $\bar{x}$, we deduce
\begin{equation}\label{eq:opt1}
 Q\bar{x} + c + \sum_{i\in \bar{\beta}}\bar{\lambda}_i\,(A_{i\,\bullet\,})^T  = 0.
\end{equation}
If there exists another isolated local minimizer $\wh{x}\in \mathbb{R}^n$ with a multiplier $\wh{\lambda}$ such that
$\mbox{supp}(\wh{\lambda}) = \bar{\beta}$,  we also have
\begin{equation}\label{eq:opt2}
Q\wh{x} + c+ \sum_{i\in \bar{\beta}}\wh{\lambda}_i\,(A_{i\,\bullet\,})^T  = 0.
\end{equation}
Multiplying both sides of \eqref{eq:opt1} by $(\bar{x} - \wh{x})^T$ and those of \eqref{eq:opt2} by $(\wh{x} - \bar{x})^T$, and by noting that
${A}_{i\,\bullet\,}\bar{x} ={A}_{i\,\bullet\,}\wh{x} = b_i$,
we may derive
\[
(\wh{x} - \bar{x})^TQ(\wh{x} - \bar{x}) = 0.
\]
Denote $\mathcal{I}(\bar{x})\,\triangleq\,\{i\mid A_{i\,\bullet\,} \bar{x} = b_i\}$.
We may write the critical cone of the problem \eqref{opt:qp} at $\bar{x}$ based on the multiplier $\bar{\lambda}$ as
\[
\mathcal{C}(\bar{x})\,\triangleq \,
\left\{\, v\in \mathbb{R}^n  \mid
{A}_{i\,\bullet\,} v \leq 0, \;\, \forall \; i\in \mathcal{I}(\bar{x});\epc
  {A}_{i\,\bullet\,}v = 0, \;\, \forall \; i \in \bar{\beta}\,
 \right\}.
\]
Since for any $i\in \mathcal{I}(\bar{x})$, ${A}_{i\,\bullet\,}\wh{x} \leq b_i$ and ${A}_{i\,\bullet\,}\bar{x} = b_i$, and for any $i\in \bar{\beta}$,
${A}_{i\,\bullet\,}\wh{x} = {A}_{i\,\bullet\,}\bar{x} = b_i$, we deduce that
$0\neq \wh{x} - \bar{x}\in \mathcal{C}(\bar{x})$.  This leads to a contradiction with the
 second order sufficient condition at the isolated local minimizer $\bar{x}$.  Therefore, the set of all isolated local minima of \eqref{opt:qp} is finite
 because the family ${\cal F}$ is finite.
\hfill $\Box$

As mentioned before, part (a) the following corollary is immediate. Part (b) is a result recently proved in \cite{CuiPang18}.
Note that a directional stationary value is derived from
a first-order directional stationary point that is not necessarily a local minimizer of the problem.

\begin{corollary} \label{co:finiteness PLQ} \rm
Let $f$ be a PLQ function on $\mathbb{R}^n$ and $X$ be a polyhedral set.  The following two statements hold for
the program (\ref{eq:basic opt problem}):
\\[0.05in]
$\bullet $ it has a finite number of isolated (strict, strong) local minima;
\\[0.05in]
$\bullet $ it has a finite number of directional stationary values.  \hfill $\Box$
\end{corollary}

The two conclusions in Corollary \ref{co:finiteness PLQ} can be obtained by invoking a very powerful finite-connected-component property of
globally subanalytic sets  \cite{BolteDaniilidisLewis06}.  This can be argued by first verifying, with a small effort, that the set of stationary 
solutions of a PLQ program
is globally subanalytic.  By the said property, it follows readily that the set of isolated stationary points must be finite.  To advance this finiteness
result to the same for strong, strict, and isolated local minima is then immediate due to their equivalence and the fact that they must be isolated
stationary points for PLQ problems.  Our proof in Proposition~\ref{pr:finite strong} is elementary,
however, and highlights one consequence of the necessity of the second-order sufficient conditions for such minima.  It is known \cite[Lemma~1.1]{FischerMarshall13}
that a PQ function on a semialgebraic set is a semialgebraic function; thus it follows from \cite{BolteDaniilidisLewis06} that a linearly constrained PQ
program  (\ref{eq:basic opt problem}) with the objective $f$ being a PQ function defined on the entire space must have finitely many isolated stationary points.  
However, it is not clear if this is sufficient to yield that this problem must have finitely many strong, strict, or isolated local minima.  
Again, the PLQ property seems needed for the latter finiteness result to hold.  

\subsection{Testing copositivity: One negative eigenvalue} \label{subsec:one negative}

Theorems~\ref{th:qp SONC} and \ref{th:qp SOSC} have shown that the (strong) local minimality of a PLQ program
can be verified via the matrix (strict) copositivity on the pieces.  The latter property can be posed
in the context of the following homogeneous quadratic program:
\begin{equation} \label{eq:homogeneous QP}
\displaystyle{
\operatornamewithlimits{\mbox{minimize}}_{v \in {\cal C}}
} \ \thalf \, v^TQv,
\end{equation}
where ${\cal C}$ is a polyhedral cone in $\mathbb{R}^n$ and $Q$ is a symmetric matrix.  The copositive of $Q$ on ${\cal C}$
then becomes the question of where the optimal objective value of (\ref{eq:homogeneous QP}) is equal to zero or unbounded below.
Since the classic work \cite{Sahni74,garey79}, it is known that a general indefinite quadratic program is NP-complete \cite{Vavasis90}.
This problem remains NP-hard even when the matrix $Q$ has only a single negative eigenvalue \cite{Pardalos91}.  In the transformations
provided in these references, the right-hand side constant in the constraint and the linear term vector in the objective are
both nonzero; this is in contrast to the problem (\ref{eq:homogeneous QP}) above which is a homogeneous problem.  Interestingly, the
homogeneity of the problem turns the hardness result in the latter reference into a computationally tractable problem.
In this subsection, we discuss the problem (\ref{eq:homogeneous QP}) when $Q$ has only one negative eigenvalue and show that
the resolution of the unboundedness of this QP
can be accomplished by solving 2 convex quadratic programs, provided that an eigen-decomposition of $Q$ is available.  As the
second-order stationarity condition of a QP, this case is related to the quasi-convexity of the
objective function; this connection is due to the known fact in generalized convexity that
the Hessian matrix of a twice differentiable quasi-convex function has only one negative eigenvalue \cite{Crouzeix99}.  In spite of this
known fact, the derivation below, although easy, does not seem to exist in the vast literature on this subject.

We begin by factoring the matrix $Q = P^{-1} D P$ where $P$ is
an orthogonal matrix whose columns are the normalized eigenvectors of $Q$, and $D$ is a diagonal matrix of eigenvalues which we denote
$\sigma_i$, for $i = 1, \cdots, n$.  Without loss of generality, we assume $\displaystyle{
\min_{1 \leq i \leq n - 1}
} \, \sigma_i \geq 0 > \sigma_n$.  With the substitution of variables $x = Pv$, the QP (\ref{eq:homogeneous QP}) is equivalent to:
\begin{equation} \label{eq:equivalent AVQP}
\begin{array}{ll}
\displaystyle{
\operatornamewithlimits{\mbox{minimize}}_{x, \, v}
} & \thalf \, x^TDx \, = \,  \underbrace{\thalf \,\displaystyle{
\sum_{i=1}^{n-1}
} \, \sigma_i \, x_i^2}_{\mbox{\small (+)ve sum of squares}} - \thalf \, | \, \sigma_n \, | \, x_n^2 \\ [0.35in]
\mbox{subject to} & \underbrace{v \, \in \, {\cal C}
\epc \mbox{and} \epc x \, = \, Pv}_{\mbox{\small remains a polyhedral cone in $(x,v)$-space}}.
\end{array} \end{equation}
Consider two related convex quadratic programs: 
\begin{equation} \label{eq:equivalent AVQP-II}
\left\{ \begin{array}{ll}
\displaystyle{
\operatornamewithlimits{\mbox{minimize}}_{y, \, \wh{v}}
} & \thalf \, \displaystyle{
\sum_{i=1}^{n-1}
} \, \sigma_i \, y_i^2 - \thalf \, | \, \sigma_n \, | \\ [0.2in]
\mbox{subject to} & \wh{v} \, \in \, {\cal C},
\epc y \, = \, P \, \wh{v}, \epc \mbox{and} \epc y_n \, = \, 1;
\end{array} \right. \end{equation}
\begin{equation} \label{eq:equivalent AVQP-III}
\left\{ \begin{array}{ll}
\displaystyle{
\operatornamewithlimits{\mbox{minimize}}_{y, \, \wh{v}}
} & \thalf \, \displaystyle{
\sum_{i=1}^{n-1}
} \, \sigma_i \, y_i^2 - \thalf \, | \, \sigma_n \, | \\ [0.2in]
\mbox{subject to} & \wh{v} \, \in \, {\cal C},
\epc y \, = \, P \, \wh{v}, \epc \mbox{and} \epc y_n \, = \, -1.
\end{array} \right. \end{equation}
Notice that the objective functions of (\ref{eq:equivalent AVQP-II}) and (\ref{eq:equivalent AVQP-III})
are bounded below on the respective feasible sets, which may be empty.  Hence if either one of these programs is feasible,
then it must attain an optimal solution.
We have the following result that connects the nonconvex QP (\ref{eq:equivalent AVQP}) with these two convex QPs
(\ref{eq:equivalent AVQP-II}) and (\ref{eq:equivalent AVQP-III}).

\begin{proposition} \label{pr:AVQP} \rm
Suppose that $Q$ has only one negative eigenvalue.  The non-convex (\ref{eq:homogeneous QP}) is unbounded below
if and only if either
(\ref{eq:equivalent AVQP-II}) or (\ref{eq:equivalent AVQP-III}) is feasible and attains a negative optimal objective value.
\end{proposition}

\noindent{\bf Proof.}  Suppose (\ref{eq:homogeneous QP}) is unbounded below.  Then there exists a feasible pair $(x,v)$ such that
the objective value of the QP (\ref{eq:equivalent AVQP}) is negative.  Clearly $x_n \neq 0$.
If $x_n > 0$, then $( y, \wh{v} ) \triangleq \displaystyle{
\frac{1}{x_n}
} \, ( x,v )$ is feasible to (\ref{eq:equivalent AVQP-II}) and its minimum objective value must be attained and is negative.
Similarly for $x_n < 0$.
Conversely, if either (\ref{eq:equivalent AVQP-II}) and (\ref{eq:equivalent AVQP-III}) has a negative optimum
objective value, then the corresponding optimal solution provides a feasible solution to (\ref{eq:equivalent AVQP}) with
a negative objective value.  Scaling this solution shows that (\ref{eq:equivalent AVQP}) is unbounded below.
\hfill $\Box$

\noindent{\bf Discussion.}  Admittedly, the materials in this subsection are so easy that we find it surprising not being able to locate the procedure
in the existing literature.  The closest result is in the reference \cite{Jargalsaikhan13} where the author considered the ``standard'' copositivity
problem on the nonnegative orthant and derived two convex ``quadratic programs''
over the second-order (Lorentz) cone whose solutions would resolve the copositivity decision problem.  In theory, the test in the reference
can be applied to any polyhedral cone provided that the generators of the cone are known, or possibly by a direct extension without invoking
such generators; neither approach is discussed, however.
Moreover, the former procedure would not be practically viable except
for special polyhedral cones.  In contrast, our procedure requires solving two standard convex quadratic programs with linear constraints
and does not require any information about the generators of the cone.  Furthermore, the procedure in Subsection~\ref{subsec:one negative}
can be extended to matrices
with exactly two negative eigenvalues, by the use of parametric convex quadratic programming \cite{CottlePangStone92} via its linear
complementarity formulation.  Nevertheless the complexity of such a parametric scheme is expected to be exponential as suggested by
the case of parametric linear programming \cite{Murty80}.  This is significantly different from the case of just one negative eigenvalue
that can be resolved by solving 2 convex quadratic programs, subject to the eigen-decomposition of the matrix in the quadratic form.
At this time, it appears that there is no practically efficient procedure for testing matrix-copositivity, thus the second-order necessary
and sufficient conditions for PLQ programs, except via the general method
of copositive programming; further research is needed.

\section{Statistical Optimization Problems} \label{subsec:statistical optimization}

In this section, we present some
modern statistical estimation problems defined by various
estimation, loss, and sparsity functions and ascertain that the objective function of the resulting optimization problem is PLQ.
This leads to the special class of problems (\ref{eq:special SC1}) that we will study in greater detail in the remaining sections.
For more details of this unified treatment of the statistical estimation problems, see \cite{CuiPangSen18}.

\noindent\underline{Piecewise affine statistical model}.  Extending the traditional
linear statistical estimation model, a piecewise affine model has recently been proposed in \cite{HahnBanergjeeSen16}
and algorithms for solving the model have been developed in \cite{CuiPangSen18}:
\begin{equation} \label{eq:PL model}
y \, = \, m(x;\Theta) + \mbox{ error }, \; \mbox{where} \;
m(x;\Theta) \, = \, \displaystyle{
\max_{1 \leq i \leq k_1}
} \, \left( \, ( \, a^i \, )^T x + \alpha_i \, \right) - \displaystyle{
\max_{1 \leq i \leq k_2}
} \, \left( \, ( \, b^i \, )^T x + \beta_i \, \right),
\end{equation}
for some positive integers $k_1$ and $k_2$.
The parameters to be estimated are contained in the tuple
$\Theta \, \triangleq \, \left\{ \, \left( \, a^i,\alpha_i \, \right)_{i=1}^{k_1}, \,
\left( \, b^i,\beta_i \, \right)_{i=1}^{k_2} \, \right\} \, \in \, \mathbb{R}^{(k_1+k_2)(d+1)}$ where each pair
$\left( \, a^i,\alpha_i \, \right)$ and $\left( \, b^i,\beta_i \, \right)$ are of dimension $d+1$.
The PA model (\ref{eq:PL model}) includes as a special case the training of 1-layer neural network by a piecewise affine activation function
\cite{NairHinton2010,GlorotBordesBengio2011} that corresponds to the following statistical model: with the vector $w$ and scalar $\alpha$
being the unknown coefficients:
\[
y \, = \, \sigma( w^Tx + \alpha ) + \mbox{ error }
\]
where $\sigma$ is a univariate piecewise affine function such as the {\sl rectified linear unit} (ReLU) which is simply the plus-function.

\noindent\underline{Loss functions}.  Deviating from the least-squares and other differentiable loss functions,
the following loss function may not be twice differentiable or convex.
\\[0.05in]
$\bullet $ The Huber loss: for some truncation scalar $K > 0$,
\[
\ell_K^{\, {\rm H}}(t) \, \triangleq \, \left\{ \begin{array}{ll}
t^2 & \mbox{if $| \, t \, | \, \leq \, K$} \\ [5pt]
K^2 + 2 \, K \, \left[ \, | \, t \, | - K \, \right] & \mbox{if $| \, t \, | \, \geq \, K$}.
\end{array} \right.
\]
The first derivative of this function is piecewise affine:
\[
( \, \ell_K^{\, {\rm H}} \, )^{\, \prime}(t) \, \triangleq \, \left\{ \begin{array}{ll}
2 \, t & \;\mbox{if $| \, t \, | \, \leq \, K$} \\ [5pt]
2 \, K \,\mbox{ sign}(t) & \; \mbox{if $| \, t \, | \, \geq \, K$}
\end{array} \right. \, = \, 2 \, \left[ \, \max\left( \, 0, \, -K - t \, \right) - \max\left( \, -t, \, -K \, \right) \, \right].
\]
This function $\ell_K^{\, {\rm H}}$ is convex, C$^{\, 1}$, and PLQ.
\\[0.05in]
$\bullet $ A loss function with margin: for some $\varepsilon > 0$,
\[
\ell(t) \, \triangleq \, \max\left( \, | \, t \, | - \varepsilon, \, 0 \, \right),
\]
employed in support vector machines with soft margins.  This function is convex and PA.
\\[0.05in]
$\bullet $ A truncated hinge loss function for binary classification \cite{WuLiu07,WuLiu15}: for some
scalar $s \leq 0$,
\[ 
\ell(t) \, \triangleq \, \max\left( \, 1 - t, \, 0 \, \right) - \max\left( \, s - t, \, 0 \, \right) \\ [5pt]
\, = \, \left\{ \begin{array}{ll}
0 & \mbox{if $t \, \geq \, 1$} \\ [5pt]
1 - t & \mbox{if $s \, \leq \, t \, \leq \, 1$} \\ [5pt]
1 - s & \mbox{if $t \, \leq \, s$}.
\end{array} \right.
\]
This function is neither convex (when $s < 0$) nor differentiable, but is piecewise affine.
%
%
%

\noindent\underline{Sparsity functions}.  As classified in \cite{AhnPangXin17}, these functions are of two kinds: exact
and surrogate.  The exact sparsity functions have the property that their zeros coincide with the $K$-sparse vectors
for some positive integer $K$; i.e., vectors with no more than $K$ nonzero components.  In contrast, the surrogate
sparsity functions are formed from univariate approximation of the discontinuous step function $| \, t \, |_0$.
A prominent exact sparsity function is
\[
P_{[K]}(w) \, \triangleq \,\displaystyle{
\sum_{i=1}^m
} \, | \, w_i \, | - \displaystyle{
\sum_{k=1}^K
} \, | \, w_{[k]} \, |,
\]
where $| \, w_{[k]} \, |$ is the $k$th largest of
the absolute values of the components of the $m$-vector $w$ arranged in non-increasing order: $\displaystyle{
\max_{1 \leq i \leq m}
} \, | w_i | \, \triangleq \, | w_{[1]} | \geq | w_{[2]} | \geq \cdots \geq | w_{[m]} | \triangleq \displaystyle{
\min_{1 \leq i \leq m}
} \, | w_i |$,
which is piecewise linear, non-separable in its arguments, and of the form $\| \, w \, \|_1 - h(w)$,
where $h$ is convex piecewise linear.

Unlike the above exact sparsity function, the surrogate sparsity functions are separable and can be written as
$P(w) = \displaystyle{
\sum_{i=1}^n
} \, p_i(w_i)$, where each
$p_i(t) = \alpha_i \, | \, t\, | - h_i(t)$ for some scalars $\alpha_i > 0$ with
$h_i$ being a convex function that is either differentiable with a piecewise affine derivative or
is itself a piecewise affine function.  Examples of these
functions include the {\sc SCAD} \cite{FanLi01} and {\sc MCP} \cite{Zhang10} functions, both of which
are univariate C$^{\, 1}$ PLQ; see the cited references for their expressions.

An example of a PA surrogate sparsity function is the capped (or truncated) $\ell_1$ function
given by $p_{\tau}(t) \triangleq \min\left( \, 1, \, \displaystyle{
\frac{| \, t \, |}{\tau}
} \, \right) = \displaystyle{
\frac{| \, t \, |}{\tau}
} - \max\left( \, \displaystyle{
\frac{| \, t \, |}{\tau}
} - 1, \, 0 \, \right)$ for some positive scalar $\tau > 0$.

{\bf Composite objectives in statistical estimation}.
Using any one of the above loss functions together with the standard least-squares loss function,
we obtain the following estimation problem: 
given $N$ data points $( x^i,y_i ) \in \mathbb{R}^{d+1}$, the optimization problem is
\begin{equation} \label{eq:loss estimation}
\displaystyle{
\operatornamewithlimits{\mbox{minimize}}_{\Theta}
} \ f_N(\Theta) \, \triangleq \displaystyle{
\frac{1}{N}
} \, \displaystyle{
\sum_{i=1}^N
} \, \ell( y_i - m(x^i;\Theta) ),
\end{equation}
where the objective function $f_N$ is the composite of the function
$( t_1, \cdots, t_N ) \mapsto \displaystyle{
\frac{1}{N}
} \, \displaystyle{
\sum_{i=1}^N
} \, \ell(t_i)$ with the vector PA function
$\Theta \mapsto \left( \, y_1 - m(x^1;\Theta), \cdots, y_N - m(x^N;\Theta) \, \right)$.
With the loss function $\ell$ being PLQ and the statistical model $m(x;\bullet)$ being PA,
the composite objective function $f_N$ is PLQ.
%
An alternative optimization problem derived from the log-likelihood maximization of a one-parameter exponential family of density
functions can be formulated as:
\begin{equation} \label{eq:log-likelihood} \displaystyle{
\operatornamewithlimits{\mbox{minimize}}_{\Theta}
} \ f_N^b(\Theta) \, \triangleq \, \displaystyle{
\frac{1}{N}
} \, \displaystyle{
\sum_{i=1}^N
} \, \ \left[ y_i \, m(x^i;\Theta) + b(m(x^i;\Theta)) \, \right],
\end{equation}
where examples of the univariate convex function $b(t)$ include: the square function $t^2$, the logarithmic function $\log( 1 + e^t )$, and
the exponential function $e^t$ corresponding to a Gaussian, Bernouilli, and a Poisson random variable, respectively.
Since $f_N^b$ is the composite of a convex function with a PA function, Proposition~\ref{pr:cvx with PA} is applicable to
(\ref{eq:log-likelihood}).

When a PLQ surrogate sparsity function is added to a composite loss
function, the resulting objective remains PLQ.
To illustrate, consider the following optimization problem
for a given scalar $\gamma > 0$,
\begin{equation} \label{eq:ReLu}
\displaystyle{
\operatornamewithlimits{\mbox{minimize}}_{w; \, \alpha}
} \ \displaystyle{
\frac{1}{N}
} \, \displaystyle{
\sum_{i=1}^N
} \, \ell\left( y_i - \sigma(w^Tx^i + \alpha_i) \, \right) + \gamma \, \displaystyle{
\sum_{i=1}^m
} \, \left[ \, \underbrace{\alpha_i \, | \, w_i \, | - h_i(w_i)}_{\mbox{surrogate sparsity function}} \, \right],
\end{equation}
where $\sigma$ is a univariate piecewise affine activation function, and each $h_i(w_i)$ is a univariate convex PLQ function.
In this case, the objective function is the sum of a weighted $\ell_1$-norm plus the function below:
\begin{equation} \label{eq:PA regression with sparsity}
( \, w, \alpha \, ) \, \mapsto \, \displaystyle{
\frac{1}{N}
} \, \displaystyle{
\sum_{i=1}^N
} \, \ell\left( y_i - \sigma(w^Tx^i + \alpha_i) \, \right) - \gamma \, \displaystyle{
\sum_{i=1}^m
} \, h_i(w_i),
\end{equation}
which is the composite of
the separable function $( t_1, \cdots, t_N, v_1, \cdots, v_m) \mapsto \displaystyle{
\frac{1}{N}
} \, \displaystyle{
\sum_{i=1}^N
} \, \ell(t_i) - \gamma \, \displaystyle{
\sum_{i=1}^m
} \, h_i(v_i)$ with the PA function:
\[
( \, w, \alpha \, ) \, \mapsto \, \left( \, y_i - \sigma(w^Tx^1 +\alpha_1), \, \cdots, \, y_N - \sigma(w^Tx^N +\alpha_N), \,
w_1, \cdots, \cdots, w_m \, \right).
\]

\section{A Class of Unconstrained Composite Programs}  \label{sec:a class of composite}

Motivated by the statistical estimation problem (\ref{eq:loss estimation}) augmented by a sparsity function
such as (\ref{eq:ReLu}),
we consider in this section a
class of unconstrained composite optimization problems and study their second-order optimality conditions:
\begin{equation} \label{eq:special SC1}
\displaystyle{
\operatornamewithlimits{\mbox{minimize}}_{w \in \mathbb{R}^n}
} \ \theta(w) \, \triangleq \, \underbrace{f(\Phi(w))}_{\mbox{denoted $\varphi(w)$}} + \, \displaystyle{
\sum_{i=1}^n
} \, \alpha_i \, | \, w_i \, | 
\end{equation}
where $f$ is a C$^{\, 1}$ PLQ function defined on $\mathbb{R}^m$ for some positive integer $m$;
$\Phi$ is a $m$-dimensional vector PA function; and each $\alpha_i$ is a nonnegative scalar.
For simplicity, we assume that
the gradient $F(z) \triangleq \nabla f(z)$ is piecewise affine with affine pieces $\left\{ A^{\, j}z + p^{\, j} \right\}_{j=1}^J$ for some positive integer $J$,
matrices $A^{\, j} \in \mathbb{R}^{m \times m}$, and vectors $e^{\, j} \in \mathbb{R}^m$; 
we further assume that $\Phi$ is PA with affine pieces
$\left\{ B^{\, k}w + q^{\, k} \right\}_{k=1}^K$ for some positive integer $K$.  This setting allows us to focus on the nondifferentiable
piecewise function $\Phi$ and the absolute-value function.
For a given $\bar{w}$, write $\bar{z} \triangleq \Phi(\bar{w})$.  Let
\[
{\cal P}_F(\bar{z}) \, \triangleq \, \left\{ \, j \, \mid \, F(\bar{z}) = A^{\, j} \bar{z} + p^{\, j} \, \right\}
\epc \mbox{and} \epc
{\cal P}_\Phi(\bar{w}) \, \triangleq \, \left\{ \, k \, \mid \, \Phi(\bar{w}) = B^{\, k} \bar{w} + q^{\, k} \, \right\}
\]
denote the active pieces of $F$ and $\Phi$ at $\bar{z}$, and $\bar{w}$, respectively.
By Proposition~\ref{pr:SC1 is dd2}, we have,
for every $v \in \mathbb{R}^n$,
\[ \begin{array}{llll}
& \theta^{\, \prime}(\bar{w};v) & = & F(\bar{z})^T \Phi^{\, \prime}(\bar{w};v) + \displaystyle{
\sum_{i \, \mid \, \bar{w}_i = 0}
} \, \alpha_i \, | \, v_i \, | + \displaystyle{
\sum_{i \, \mid \, \bar{w}_i \neq 0}
} \, \alpha_i \, v_i \,  \mbox{sign}(\bar{w}_i)  \\ [0.2in]
\mbox{and} &
\theta^{\, (2)}(\bar{w};v) & = & f^{(2)}(\bar{z};\Phi^{\, \prime}(\bar{w};v)) \, = \,
\Phi^{\, \prime}(\bar{w};v)^TF^{\, \prime}(\bar{z};\Phi^{\, \prime}(\bar{w};v)).
\end{array}
\]
Since $F^{\prime}(\bar{z};\bullet)$ and $\Phi^{\, \prime}(\bar{w};\bullet)$ are PL functions, it follows,
by \cite[Lemma~4.6.1]{FacchineiPang03}, that for every $u \in \mathbb{R}^m$ and $v \in \mathbb{R}^n$, there exist subsets $\wh{\cal P}_F(\bar{z};u)$ and
$\wh{\cal P}_{\Phi}(\bar{w};v)$ of ${\cal P}_F(\bar{z})$ and ${\cal P}_{\Phi}(\bar{w})$, respectively, such that
\[ \begin{array}{llll}
& F^{\, \prime}(\bar{z};u) & = & A^{\, j}u \epc \forall \, j \, \in \, \wh{\cal P}_F(\bar{z};u) \\ [0.15in]
\mbox{and} & \Phi^{\, \prime}(\bar{w};v) & = & B^{\, k}v \epc \forall \, k \, \in \, \wh{\cal P}_{\Phi}(\bar{w};v).
\end{array} \]
These index sets $\wh{\cal P}_F(\bar{z};u)$ and
$\wh{\cal P}_{\Phi}(\bar{w};v)$ contain the directionally active indices like ${\cal A}^{\, \prime}(x;d)$ for a general PA function.
The following result is an immediate consequence of Theorem~\ref{th:qp SOSC}, giving necessary and sufficient conditions for $\bar{w}$ to be
a (strong, isolated, strict) local minimizer of (\ref{eq:special SC1}) in terms of the second-order conditions (a), (b), and (c).

\begin{proposition} \label{pr:2nd order statistical program} \rm
Consider the following three conditions:

(a) $F(\bar{z})^T B^{\, k}v + \displaystyle{
\sum_{i \, : \, \bar{w}_i \neq 0}
} \, \alpha_i \, v_i \, \mbox{sign}(\bar{w}_i) + \displaystyle{
\sum_{i \, : \, \bar{w}_i = 0}
} \, \alpha_i \, | \, v_i \, | \geq 0$ for all $v \in \mathbb{R}^n$ and all $k \in \wh{\cal P}_{\Phi}(\bar{w};v)$;

(b) for all $v \in \mathbb{R}^n$ and all pairs $(j,k) \in \wh{\cal P}_F(\bar{z};u) \times \wh{\cal P}_{\Phi}(\bar{w};v)$
where $u = \Phi^{\, \prime}(\bar{w};v)$,
\[
\left[ F(\bar{z})^T B^{\, k}v + \displaystyle{
\sum_{i \, : \, \bar{w}_i \neq 0}
} \, \alpha_i \, v_i \, \mbox{sign}(\bar{w}_i) + \displaystyle{
\sum_{i \, : \, \bar{w}_i = 0}
} \, \alpha_i \, | \, v_i \, | \leq  0 \right] \ \Rightarrow \
v^T \left[ \, ( \, B^{\, k} \, )^T A^{\, j} B^{\, k} \, \right] \, v \, \geq \, 0;
\]
(c) for all $v \in \mathbb{R}^n$ and all pairs $(j,k) \in \wh{\cal P}_F(\bar{z};u) \times \wh{\cal P}_{\Phi}(\bar{w};v)$
where $u = \Phi^{\, \prime}(\bar{w};v)$,
\begin{align*}
&\left[ \, v \, \neq \, 0 \mbox{ and } F(\bar{z})^T B^{\, k}v + \displaystyle{
	\sum_{i \, : \, \bar{w}_i \neq 0}
} \, \alpha_i \, v_i \, \mbox{sign}(\bar{w}_i) + \displaystyle{
\sum_{i \, : \, \bar{w}_i = 0}
} \, \alpha_i \, | \, v_i \, | \, \leq \, 0 \, \right] \\
& \Rightarrow
v^T \left[ \, ( \, B^{\, k} \, )^T A^{\, j} B^{\, k} \, \right] \, v \, > \, 0.
\end{align*}
It holds that
\\[0.05in]
$\bullet $ conditions (a) and (b) combined are necessary and sufficient for $\bar{w}$ to be a local minimizer of (\ref{eq:special SC1});
\\[0.05in]
$\bullet $ conditions (a) and (c) combined are necessary and sufficient for $\bar{w}$ to be a strong (equivalently, strict or isolated)
local minimizer of (\ref{eq:special SC1});
\\[0.05in]
$\bullet $ the number of strong (strict, or isolated) local minimizers is finite;
\\[0.05in]
$\bullet $ the number of directional stationary values is finite.  \hfill $\Box$
\end{proposition}

Unlike the sets ${\cal P}_F(\bar{z})$ and ${\cal P}_\Phi(\bar{w})$ which are completely determined, respectively, by the vectors
$\bar{z}$ and $\bar{w}$ alone, elements of the sets $\wh{\cal P}_F(\bar{z};u)$ and $\wh{\cal P}_\Phi(\bar{w};v)$ cannot be totally identified based
only on the pairs $(\bar{z},u)$ and $(\bar{w},v)$, respectively.  Indeed, $\wh{\cal P}_F(\bar{z};v)$ consists of all indices $j \in {\cal P}_F(\bar{z})$
such that $F(\bar{z} + \tau_{\nu} u) = A^{\, j} ( \bar{z} + \tau_{\nu} u ) + e^{\, j}$ for a sequence of positive scalars $\{ \tau_{\nu} \} \downarrow 0$.
A similar description applies to the elements in $\wh{\cal P}_{\Phi}(\bar{w};v)$.
Thus if either ${\cal P}_F(\bar{z})$ or ${\cal P}_{\Phi}(\bar{w})$ is not a singleton, the verification of the second or third condition
in Proposition~\ref{pr:2nd order statistical program} does not appear to be easy without enumeratively checking all pairs of indices in these
index sets.  This is the combinatorial aspect of the non-smoothness of the composite function $\varphi(w)$.

\section{A Homogeneous Singly Absolute-Value Constrained QP}

Consider the simplified situation of (\ref{eq:special SC1}) where both
${\cal P}_F(\bar{z})$ and ${\cal P}_{\Phi}(\bar{w})$ are singletons.  This motivates the investigation of an indefinite quadratic
optimization problem (\ref{eq:AVQP}) with a single absolute-value constraint that aims to address the two second-order conditions
in (b) and (c) in Proposition~\ref{pr:2nd order statistical program}.
We show that the resolution of the problem (\ref{eq:AVQP}) is equivalent to testing the copositivity of a certain matrix
on a nonnegative orthant, and thus is in general NP-hard \cite{MurtyKabadi87,Vavasis90}.


Let $Q \in \mathbb{R}^{n \times n}$ be a symmetric indefinite matrix, $b \in \mathbb{R}^n$ be
arbitrary, and $\alpha \in \mathbb{R}^n$ be a nonnegative, nonzero vector.  Consider the quadratic program (QP) with a homogeneous objective:
\begin{equation} \label{eq:AVQP}
\displaystyle{
\operatornamewithlimits{\mbox{minimize}}_{v \in \mathbb{R}^n}
} \ \thalf \, v^TQv \epc \mbox{subject to} \epc  b^Tv + \displaystyle{
\sum_{i=1}^n
} \, \alpha_i \, | \, v_i \, | \, \leq \, 0,
\end{equation}
where the constraint is such that the reverse inequality holds for all vectors $v \in \mathbb{R}^n$; thus $b$ and $\alpha$ satisfy:
\[
\underbrace{\left[ \, b_i \, v_i + \alpha_i \, | \, v_i \, | \, \geq \, 0, \epc \forall \, v_i \, \in \, \mathbb{R} \, \right] \epc \forall \, i \, = \, 1, \cdots, n}_{\mbox{
this is the inequality in part (a) of Proposition~\ref{pr:2nd order statistical program} for one index $k$}},
\]
which is equivalent to $| b_i | \leq \alpha_i$ for all $i = 1, \cdots, n$.  Based on this observation, we can derive the following
lemma which shows in particular that the constraint set of (\ref{eq:AVQP}) is the Cartesian product of four types of 1-dimensional rays:
$\{ 0 \}$ (a degenerate ray); the entire real line, the nonnegative, or nonpositive real axis.

\begin{lemma} \label{lm:constraint signs} \rm
Let $b$ and $\alpha$ be $n$-vectors such that $| b_i | \leq \alpha_i$ for all $i = 1, \cdots, n$.  A vector $v \in \mathbb{R}^n$ satisfies
(\ref{eq:AVQP}) if and only if the following three conditions hold for all $i = 1, \cdots, n$,\\[0.05in]
$\bullet $ $| b_i | < \alpha_i$ implies $v_i = 0$
\\[0.05in]
$\bullet $ $| b_i | = \alpha_i > 0$ implies either $v_i  = 0$ or $\mbox{sign}(v_i) = -\mbox{sign}(b_i)$;
\\[0.05in]
$\bullet $ $b_i = \alpha_i = 0$ implies $v_i$ free.
\end{lemma}
{\bf Proof.}  We can write: 
\[ \begin{array}{lll}
b^Tv + \displaystyle{
\sum_{i=1}^n
} \, \alpha_i \, | \, v_i \, |  & = & \displaystyle{
\sum_{i \, : \, | b_i | < \alpha_i}
}   \left( \, b_i \, v_i + \alpha_i \, | \, v_i \, | \, \right) + \displaystyle{
\sum_{i \, : \, | b_i | = \alpha_i > 0}
}  \left( \, b_i \, v_i + \alpha_i \, | \, v_i \, | \, \right) + \displaystyle{
\sum_{i \, : \, b_i = \alpha_i = 0}
} \left( \, b_i \, v_i + \alpha_i \, | \, v_i \, | \, \right) \\ [0.25in]
& = & \displaystyle{
\sum_{i \, : \, | b_i | < \alpha_i}
}  \left( \, b_i \, v_i + \alpha_i \, | \, v_i \, | \, \right) + \displaystyle{
\sum_{i \, : \, | b_i | = \alpha_i > 0}
} \left( \, b_i \, v_i + | \, b_i \, v_i \, | \, \right).
\end{array} \] 
Thus, $b^Tz + \displaystyle{
\sum_{i=1}^n
} \, \alpha_i \, | \, v_i \, | \leq 0$ if and only if each of the summands on the right-hand side is equal to zero.
This readily yields the desired equivalence.  \hfill $\Box$

Before proceeding further, we mention that although this section has focused on the QP (\ref{eq:AVQP}) with one single convex
absolute-value constraint, it is easy to generalize the analysis to arbitrary linear constraints.   The end result is that
we can obtain similar characterizations of the second-order conditions for PLQ programs in terms of certain matrix-copositivity
properties of Schur complements on the nonnegative orthant.

Under the assumption that
$| b_i | \leq \alpha_i$ for all $i = 1, \cdots, n$, the problem (\ref{eq:AVQP}) is thus equivalent to
\begin{equation} \label{eq:equivalent AVQP-0}
\displaystyle{
\operatornamewithlimits{\mbox{minimize}}_{z \in \mathbb{R}^n}
} \ \thalf \, v^TQv \epc \mbox{subject to}  \left\{ \begin{array}{lll}
v_i \, = \, 0 & \mbox{if $| \, b_i \, | < \alpha_i$} \\ [5pt]
v_i \, \geq \, 0 & \mbox{if $-b_i = \alpha_i > 0$} & \mbox{index set denoted ${\cal I}_+$} \\ [5pt]
v_i \, \leq \, 0 & \mbox{if $b_i = \alpha_i > 0$} & \mbox{index set denoted ${\cal I}_-$} \\ [5pt]
v_i \mbox{ free} & \mbox{if $b_i = \alpha_i = 0$} & \mbox{index set denoted ${\cal I}_f$}.
\end{array} \right.
\end{equation}

This homogeneous program is either unbounded below or has a zero optimum objective value.  The latter happens
if and only if the matrix 
\begin{equation} \label{eq:partitioned copositive-I}
\left[ \begin{array}{cccc}
 Q_{{\cal I}_+ {\cal I}_+} & -Q_{{\cal I}_+ {\cal I}_-} & | &  Q_{{\cal I}_+ {\cal I}_f} \\ [5pt]
-Q_{{\cal I}_- {\cal I}_+} &  Q_{{\cal I}_- {\cal I}_-} & | & -Q_{{\cal I}_- {\cal I}_f} \\
\rule{0.5in}{0.01in} & \rule{0.5in}{0.01in} & | & \rule{0.5in}{0.01in} \\ [5pt]
 Q_{{\cal I}_f {\cal I}_+} & -Q_{{\cal I}_f {\cal I}_-} & | & Q_{{\cal I}_f {\cal I}_f}
\end{array} \right]
\end{equation} 
is copositive on the ``mixed cone'' $\mathbb{R}_+^{| {\cal I}_{\neq 0} |} \times \mathbb{R}^{| {\cal I}_f |}$,
where ${\cal I}_{\neq 0} \triangleq {\cal I}_+ \cup {\cal I}_-$.  In what follows, we perform matrix operations
to remove the subspace $\mathbb{R}^{| {\cal I}_f |}$ and
convert this copositivity condition on the mixed cone into the copositivity of a matrix of reduced order on the
nonnegative orthant $\mathbb{R}_+^{| {\cal I}_{\neq 0} |}$.  We begin by noting that a necessary condition for
the copositivity of the matrix (\ref{eq:partitioned copositive-I})
on the mixed cone is that the submatrix $Q_{{\cal I}_f {\cal I}_f}$ is positive semidefinite.  As such, there exist
an orthogonal matrix $P_{{\cal I}_f {\cal I}_f}$ of order $| {\cal I}_f |$ of normalized eigenvectors of $Q_{{\cal I}_f {\cal I}_f} $
and a diagonal matrix of $\Xi_{{\cal I}_f}$ with nonnegative diagonals such that
$\left[ \, P_{{\cal I}_f {\cal I}_f} \, \right]^T Q_{{\cal I}_f {\cal I}_f} \left[ \, P_{{\cal I}_f {\cal I}_f} \, \right] =  \Xi_{{\cal I}_f}$.
It is not difficult to show that the matrix (\ref{eq:partitioned copositive-I}) is copositive on
$\mathbb{R}_+^{| {\cal I}_{\neq 0} |} \times \mathbb{R}^{| {\cal I}_f |}$ if and only if $Q_{{\cal I}_f {\cal I}_f}$ is positive semidefinite
and the matrix 
\begin{equation} \label{eq:partitioned copositive-II}
\left[ \begin{array}{cccc}
 Q_{{\cal I}_+ {\cal I}_+} & -Q_{{\cal I}_+ {\cal I}_-} & | &  Q_{{\cal I}_+ {\cal I}_f} \, P_{{\cal I}_f {\cal I}_f} \\ [5pt]
-Q_{{\cal I}_- {\cal I}_+} &  Q_{{\cal I}_- {\cal I}_-} & | & -Q_{{\cal I}_- {\cal I}_f} \, P_{{\cal I}_f {\cal I}_f} \\
\rule{0.7in}{0.01in} & \rule{0.7in}{0.01in} & | & \rule{0.7in}{0.01in} \\ [5pt]
\left[ \, P_{{\cal I}_f {\cal I}_f} \, \right]^T Q_{{\cal I}_f {\cal I}_+} & -\left[ \, P_{{\cal I}_f {\cal I}_f} \, \right]^T Q_{{\cal I}_f {\cal I}_-} & | & \Xi_{{\cal I}_f}
\end{array} \right]
\end{equation} 
is copositive on the same cone.  We may partition the index set ${\cal I}_f$ into the union of two complementary index subsets
${\cal I}_f^+$ and ${\cal I}_f^0$ such that $\Xi_{{\cal I}_f} = \left[ \begin{array}{cc}
\Xi_{{\cal I}_f^+} & 0 \\
0 & 0
\end{array} \right]$ where $\Xi_{{\cal I}_f^+}$ is a diagonal matrix with positive diagonals.  These
preparatory manipulations lead to the following reduction result for the quadratic form $v^TQv$ to be
nonnegative on the feasible set of (\ref{eq:AVQP}) under the given stipulation of the coefficients
$b_i$ and $\alpha_i$.

\begin{proposition} \label{pr:nonnegative quadratic form} \rm
Suppose $| b_i | \leq \alpha_i$ for all $i = 1, \cdots, n$.  A necessary and sufficient condition for the
quadratic program (\ref{eq:AVQP}) to have a zero optimum objective value is for the three conditions below
to hold:

\noindent$\bullet $ the principle submatrix $Q_{{\cal I}_f {\cal I}_f}$ is positive semidefinite with eigen-decomposition \\
$\left[ \, P_{{\cal I}_f {\cal I}_f} \, \right]^T Q_{{\cal I}_f {\cal I}_f} \left[ \, P_{{\cal I}_f {\cal I}_f} \, \right] =  \Xi_{{\cal I}_f}$;

\noindent$\bullet $ $\left[ \begin{array}{c}
Q_{{\cal I}_+ {\cal I}_f} \\
Q_{{\cal I}_- {\cal I}_f}
\end{array} \right] P_{{\cal I}_f {\cal I}_f^0} = 0$;

\noindent$\bullet $ the Schur complement 
\[
\left[ 
\begin{array}{cc}
 Q_{{\cal I}_+ {\cal I}_+} & -Q_{{\cal I}_+ {\cal I}_-} \\ [5pt]
-Q_{{\cal I}_- {\cal I}_+} &  Q_{{\cal I}_- {\cal I}_-}
\end{array} 
\right] - \left[ \begin{array}{c}
Q_{{\cal I}_+ {\cal I}_f} \, P_{{\cal I}_f {\cal I}_f^+} \\ [5pt]
-Q_{{\cal I}_- {\cal I}_f} \, P_{{\cal I}_f {\cal I}_f^+}
\end{array} \right] \, 
\left[ \, \Xi_{{\cal I}_f^+} \, \right]^{-1}
\left[ \begin{array}{c}
Q_{{\cal I}_+ {\cal I}_f} \, P_{{\cal I}_f {\cal I}_f^+} \\ [5pt]
-Q_{{\cal I}_- {\cal I}_f} \, P_{{\cal I}_f {\cal I}_f^+}
\end{array} \right]^T
\] 
is copositive on $\mathbb{R}_+^{| {\cal I}_{\neq 0} |}$.
\end{proposition}
{\bf Proof.}  ``Necessity''.  The matrix (\ref{eq:partitioned copositive-II}) can be written in further partitioned form: 
\begin{equation} \label{eq:partitioned copositive-II-1}
\left[ \begin{array}{ccccc}
 Q_{{\cal I}_+ {\cal I}_+} & -Q_{{\cal I}_+ {\cal I}_-} & | &  Q_{{\cal I}_+ {\cal I}_f} \, P_{{\cal I}_f {\cal I}_f^+} &
 Q_{{\cal I}_+ {\cal I}_f} \, P_{{\cal I}_f {\cal I}_f^0} \\ [5pt]
-Q_{{\cal I}_- {\cal I}_+} &  Q_{{\cal I}_- {\cal I}_-} & | & -Q_{{\cal I}_- {\cal I}_f} \, P_{{\cal I}_f {\cal I}_f^+} &
-Q_{{\cal I}_+ {\cal I}_f} \, P_{{\cal I}_f {\cal I}_f^0} \\ [5pt]
\rule{0.8in}{0.01in} & \rule{0.8in}{0.01in} & | & \rule{0.8in}{0.01in} & \rule{0.8in}{0.01in} \\ [5pt]
\left[ \, P_{{\cal I}_f {\cal I}_f^+} \, \right]^T Q_{{\cal I}_f {\cal I}_+} & -\left[ \, P_{{\cal I}_f {\cal I}_f^+} \, \right]^T Q_{{\cal I}_f {\cal I}_-}
& | & \Xi_{{\cal I}_f^+} & 0 \\ [0.1in]
\left[ \, P_{{\cal I}_f {\cal I}_f^0} \, \right]^T Q_{{\cal I}_f {\cal I}_+} & -\left[ \, P_{{\cal I}_f {\cal I}_f^0} \, \right]^T Q_{{\cal I}_f {\cal I}_-}
& | & 0 & 0
\end{array} \right] .
\end{equation} 
For the latter symmetric matrix to be copositive on the mixed cone $\mathbb{R}_+^{| {\cal I}_{\neq 0} |} \times \mathbb{R}^{| {\cal I}_f |}$, it
is necessary that $\left[ \begin{array}{c}
Q_{{\cal I}_+ {\cal I}_f} \\
Q_{{\cal I}_- {\cal I}_f}
\end{array} \right] P_{{\cal I}_f {\cal I}_f^0} = 0$.  To prove the copositivity of the Schur complement, let
$v_{{\cal I}_+}$ and $v_{{\cal I}_-}$ be arbitrary nonnegative vectors.  Let
\[
v_{{\cal I}_f^+} \, \triangleq \, -\left[ \,  \Xi_{{\cal I}_f^+} \, \right]^{-1} \left\{
\left[ \, P_{{\cal I}_f {\cal I}_f^+} \, \right]^T Q_{{\cal I}_f {\cal I}_+} v_{{\cal I}_+} -
\left[ \, P_{{\cal I}_f {\cal I}_f^+} \, \right]^T Q_{{\cal I}_f {\cal I}_-} v_{{\cal I}_-}
\right\}.
\]
We then have
\[ \begin{array}{lll}
0 & \leq & \left( \, \begin{array}{c}
v_{{\cal I}_+} \\ [5pt]
v_{{\cal I}_-} \\ [5pt]
v_{{\cal I}_f^+}
\end{array} \right)^T \left[ \begin{array}{cccc}
 Q_{{\cal I}_+ {\cal I}_+} & -Q_{{\cal I}_+ {\cal I}_-} & \hspace{-0.3cm} | &  \hspace{-0.2cm} Q_{{\cal I}_+ {\cal I}_f} \, P_{{\cal I}_f {\cal I}_f^+}  \\ [5pt]
-Q_{{\cal I}_- {\cal I}_+} &  Q_{{\cal I}_- {\cal I}_-} &\hspace{-0.3cm}  |  & \hspace{-0.2cm} -Q_{{\cal I}_- {\cal I}_f} \, P_{{\cal I}_f {\cal I}_f^+} \\ [5pt]
\rule{0.8in}{0.01in} & \rule{0.8in}{0.01in} &\hspace{-0.3cm}  |  & \hspace{-0.2cm} \rule{0.8in}{0.01in}\\ [5pt]
\left[ \, P_{{\cal I}_f {\cal I}_f^+} \, \right]^T Q_{{\cal I}_f {\cal I}_+} & -\left[ \, P_{{\cal I}_f {\cal I}_f^+} \, \right]^T Q_{{\cal I}_f {\cal I}_-}
& \hspace{-0.3cm} |  & \hspace{-0.2cm} \Xi_{{\cal I}_f^+}
\end{array} \right] \left(\begin{array}{c}
v_{{\cal I}_+} \\ [5pt]
v_{{\cal I}_-} \\ [5pt]
v_{{\cal I}_f^+}
\end{array} \right) \\ [0.6in]
& = &  \left( \, \begin{array}{c}
v_{{\cal I}_+} \\ [5pt]
v_{{\cal I}_-}
\end{array} \right)^T \left[ \begin{array}{cccc}
 Q_{{\cal I}_+ {\cal I}_+} & -Q_{{\cal I}_+ {\cal I}_-} & | &  Q_{{\cal I}_+ {\cal I}_f} \, P_{{\cal I}_f {\cal I}_f^+}  \\ [5pt]
-Q_{{\cal I}_- {\cal I}_+} &  Q_{{\cal I}_- {\cal I}_-} & | & -Q_{{\cal I}_- {\cal I}_f} \, P_{{\cal I}_f {\cal I}_f^+}
\end{array} \right]  \left( \, \begin{array}{c}
v_{{\cal I}_+} \\ [5pt]
v_{{\cal I}_-} \\ [5pt]
v_{{\cal I}_f^+}
\end{array} \right).
\end{array} \]
Substituting the definition of the vector $v_{{\cal I}_f^+}$ easily the completes the proof of the necessity of the third condition.

``Sufficiency''.  This can be proved by reversing the above arguments.  \hfill $\Box$



\end{document}